\documentclass[12pt, a4paper]{amsart}

\usepackage{latexsym}

\usepackage{amssymb}
\usepackage[centertags]{amsmath}
\usepackage{verbatim}
\usepackage{epic, eepic}

\usepackage{color}

\newtheorem{theorem}{Theorem}
\newtheorem{lemma}[theorem]{Lemma}
\newtheorem{proposition}[theorem]{Proposition}
\newtheorem{corollary}[theorem]{Corollary}
\newtheorem{remark}[theorem]{Remark}

\numberwithin{theorem}{section}

\newcommand{\ci}[1]{_{{}_{\scriptstyle #1}}}

\renewcommand{\phi}{\varphi}

\newcommand{\sign}{\operatorname{sgn}}

\newcommand{\U}{\overline U}

\newcommand{\eps}{\varepsilon}

\newcommand{\la}{\lambda}

\newcommand{\V}{\mathbb V}

\renewcommand{\epsilon}{\varepsilon}

\def\done{{1\hskip-2.5pt{\rm l}}}

\newcommand{\Ex}{\mathbb E}

\renewcommand{\le}{\leqslant}
\renewcommand{\ge}{\geqslant}

\newcommand{\R}{\mathbb R}
\newcommand{\C}{\mathbb C}
\newcommand{\Z}{\mathbb Z}

\newcommand{\cZ}{\mathcal Z}

\newcommand{\T}{\mathbb T}

\title
[Random complex zeroes: asymptotic normality revisited]
{Fluctuations in random complex zeroes:
\\ Asymptotic normality revisited }

\author{F. Nazarov}
\address{F.N.: Mathematics Department \\
University of Wisconsin-Madison \\
480 Lincoln Dr., Madison WI 53706\\
USA}

\email{nazarov@math.wisc.edu}

\author{M. Sodin}
\address{M.S.: School of Mathematical Sciences\\
Tel Aviv University\\
Tel Aviv 69978\\
Israel}

\email{sodin@post.tau.ac.il}

\thanks{F.N. is partially supported by
the National Science Foundation, DMS grant 0501067. M.S. is
partially supported by the Israel Science Foundation of the Israel
Academy of Sciences and Humanities, grant 171/07.}

\date{February 21, 2010}

\begin{document}

\begin{abstract}

The Gaussian Entire Function
\[
F(z) = \sum_{k=0}^\infty \zeta_k \frac{z^k}{\sqrt{k!}}
\]
($\zeta_0, \zeta_1, \dots $ are Gaussian i.i.d. complex
random coefficients) is distinguished by the distribution invariance
of its zero set with respect to the isometries of the complex plane.
We find close to optimal conditions on a function $h$ that yield
asymptotic normality of linear statistics of zeroes
\[
n(R, h) = \sum_{a\colon F(a)=0} h\left( \frac{a}{R}\right)
\]
when $R\to\infty$, and construct examples of functions $h$ with abnormal
fluctuations of linear statistics. We show that the fluctuations of $n(R, h)$ are
asymptotically normal when $h$ is either a $C^\alpha_0$-function with $\alpha>1$,
or a $C_0^\alpha$-function with $\alpha\le 1$ such that the variance
of $n(R, h)$ is at least $R^{-2\alpha+\epsilon}$ with some $\epsilon>0$.

These results complement our recent results from ``Clustering of correlation functions
for random complex zeroes'', where, using different methods, we prove that
the fluctuations of $n(R, h)$
are asymptotically normal when $h$ is bounded, and the variance
of $n(R, h)$ grows at least as a positive power of $R$.

\end{abstract}

\maketitle

\section{Introduction and main results}\label{section_intro}

\subsection{} Consider the Gaussian entire function (G.E.F., for short)
\[
F(z) = \sum_{k=0}^\infty \zeta_k \frac{z^k}{\sqrt{k!}}
\]
where $\zeta_n$ are standard independent Gaussian complex variables
(that is, the density of $\zeta_n$ on the complex plane is $\tfrac1{\pi} e^{-|\zeta|^2}$).
A remarkable feature of the random zero set $\cZ_F=F^{-1}\{0\}$ is its distribution
invariance with respect to the isometries of the plane. The rotation invariance of $\cZ_F$ is
obvious since the distribution of  $F$ is rotation invariant itself.
Though the distribution of $F$ is not translation invariant at all,
the translation invariance of the zero process $\cZ_F$ follows from the fact that, for every $w\in\C$, the
Gaussian processes $F(z+w)$ and $F_w(z)=F(z)e^{z\overline{w}+\frac12 |w|^2}$ have the same distributions.
The latter follows by inspection of the covariances:
\[
\Ex \left\{ F_w(z_1) \overline{F_w(z_2)} \right\} = e^{z_1\overline{z_2} + z_1\overline{w} + \overline{z_2} w +|w|^2}
= e^{(z_1+w)\overline{(z_2+w)}} = \Ex \left\{ F(z_1+w)\overline{F(z_2+w)} \right\} \,.
\]

By Calabi's rigidity \cite[Section~3]{Sodin}, \cite[Section~2.5]{HKPV}, $F(z)$
is essentially the only Gaussian functions analytic in $\C$ with the distribution of zeroes
invariant with respect to the isometries of the plane. The random zero process $\cZ_F$ together
with similar models related to other one-dimensional complex geometries was introduced
in the works of Bogomolny, Bohigas, Leboeuf~\cite{BBL}, Hannay~\cite{H} and Kostlan~\cite{Kostlan}.
After that the zero process $\cZ_F$ has been studied from different points of view.
A brief non-technical introduction can be found in~\cite{NS1}.

\subsection{Linear statistics of zeroes and their variances}
One of the ways to understand the asymptotic behaviour of the zero process $\cZ_F$ is to introduce
the random variable
\[
n(R, h) = \sum_{a\in \cZ_F} h\left( \frac{a}{R} \right)
\]
called the {\em linear statistics} of $\cZ_F$ and to study its behavior as $R\to\infty$. Here,
$h$ is a function that we always assume real-valued, not identically zero,
and belonging to $ (L^1 \cap L^2)(\R^2) $. Since
$\cZ_F$ has a translation-invariant distribution, we have
\[
\Ex \left\{ n(R, h) \right\} = CR^2 \int_{\R^2} h\,.
\]
A straightforward computation shows that $C = \tfrac1{\pi}$ (see, for instance, \cite[Part~{\rm I}]{ST}).

Denote by \[ \V (R, h) =  \Ex \bigl\{ n(R, h) - \Ex n(R, h) \bigr\}^2 \]
the variance of  $n(R, h)$, and by $ \sigma (R, h) = \sqrt{\V (R, h)}$ its standard deviation.
In \cite{FH}, Forrester and Honner found that if $h\in C^2_0$ (i.e., $h$ is a
$C^2$-function with compact support), then
\[
\V (R, h) = \frac{\zeta (3) + o(1)}{16\pi R^2}\, \|\Delta h\|_{L^2}^2\,, \qquad  R\to\infty\,,
\]
while for $h=\done_G$ (the indicator function of a bounded domain $G$ with piecewise smooth boundary)
\[
\V (R, \done_G) = \frac{\zeta(3/2)+o(1)}{8\pi^{3/2}}\, R\, L(\partial G)\,, \qquad R\to\infty\,.
\]
Here, $\zeta(\,\cdot\,)$ is Riemann's zeta-function. Later, Sodin and Tsirelson~\cite[Part~{\rm I}]{ST} and
Shiffman and Zelditch~\cite{SZ} found different approaches to asymptotic computation of the variance,
which work in more general contexts. Our first result is an exact formula for the variance
$\V (R, h) $ valid for {\em arbitrary} functions $h\in (L^1 \cap L^2)(\R^2)$.

In what follows, $X \lesssim Y$ means that $X \le C\cdot Y$ with a positive numerical constant $C$,
and $X\simeq Y$ means that $X \lesssim Y$ and $Y \lesssim X $ simultaneously. By
\[
\widehat{h}(\la) = \int_{\R^2} h(x) e^{-2\pi {\rm i}\, \langle \la, x \rangle }\, {\rm d}A(x)
\]
we denote the Fourier transform of $h$.
\begin{theorem}\label{thm_variance}
For every non-zero function $h\in (L^1 \cap L^2)(\R^2)$ and each $R>0$, one has
\[
\V (R, h) = R^2\int_{\R^2}|\widehat h(\la)|^2 M(R^{-1}\la)\,{\rm d}A(\la)\,
\]
where
\[
M(\la)=\pi^3|\la|^4\sum_{\alpha\ge 1}\frac 1{\alpha^3}e^{-\frac{\pi^2}{\alpha}|\la|^2}\,.
\]
\end{theorem}

Noticing that $M(\la)\simeq\min(|\la|^4,1)$, we get
\begin{equation}\label{eq_variance}
\V (R, h) \simeq R^{-2} \int_{|\la|\le R} |\widehat{h}(\la)|^2 |\la|^4 \, {\rm d}A(\la) + R^2
\int_{|\la|\ge R} |\widehat{h}(\la)|^2 \, {\rm d}A(\la) \,.
\end{equation}
Observe that the right-hand side of \eqref{eq_variance} interpolates the $L^2$-norm of
the function $h$ and the $L^2$-norm of its Laplacian $\Delta h$.

\subsection{Normality of fluctuations}
The question about the asymptotic normality of fluctuations of $n(R, h)$
is more delicate. For compactly supported $C^2$-functions $h$, the asymptotic normality was proven
by Sodin and Tsirelson in \cite[Part~{\rm I}]{ST}.
The proof was based on the method
of moments, and the moments were computed using the diagram technique%
\footnote{Note that using the variance estimate \eqref{eq_variance}, one readily extends
this result to the functions $h$ in the Sobolev space $W_2^2$; i.e., to the functions
$h$ such that
\[
\int_{\R^2} |\widehat{h}(\la)|^2 \left( 1+|\la|^4\right)\, {\rm d}A(\la) < \infty\,.
\]
}.
After some modification, the argument
from \cite[Part~{\rm I}]{ST} also works in the case when $h=\done_G$ where $G$ is a
bounded domain with piecewise smooth boundary. Krishnapur noticed in
his PhD Thesis \cite{Krishnapur}
that the same argument also works when $h=\done_G * \done_G$ (the convolution square).
In the latter case, the variance tends to a constant. In a recent work \cite{Tsirel},
Tsirelson found the asymptotics of the logarithm of the characteristic functional
$\Ex e^{\lambda n(R, h)}$, where $h$ is a compactly supported $C^2$-function,
when $R\to\infty$, $\lambda\to  0$ in such way that
$\lambda \log^2 R \to 0$. Among other things, this gives
a different proof of asymptotic normality of smooth linear statistics of random complex zeroes.

Let $C^\alpha_0$ be the class of compactly supported $C^\alpha$-functions.
In this paper, we explore what happens for $C^\alpha_0$-functions with $0<\alpha<2$.
We say that the linear statistics $n(R, h)$ have {\em asymptotically normal fluctuations} if
the normalized linear statistics
\[
\frac{n(R, h) - \Ex n(R, h)}{\sigma (R, h)}
\]
converge in distribution to the standard (real) Gaussian random variable as $R\to\infty$.

\begin{theorem}\label{thm_normality}
Suppose that $h\in C^\alpha_0$ with some $\alpha >0$, and that for some $\epsilon>0$
and each sufficiently big $R$, we have
\begin{equation}\label{eq_ii}
\sigma (R, h)>R^{-\alpha+\epsilon}\,.
\end{equation}
Then the linear statistics $n(R, h)$ have asymptotically normal fluctuations.
\end{theorem}

Note that by \eqref{eq_variance}, we always have $\sigma (R, h) \gtrsim c(h) R^{-1}$
with positive $ c(h) $ independent of $R$. Hence, for $\alpha>1$, condition \eqref{eq_ii}
holds automatically.

\begin{corollary}\label{cor_1}
Suppose that $h\in C^\alpha_0$ with $\alpha>1$.
Then the linear statistics
$n(R, h)$ have asymptotically normal fluctuations.
\end{corollary}

We mention that there are plenty of nice $C^\alpha_0$-functions with $\alpha\le 1$
satisfying condition~\eqref{eq_ii} of Theorem~\ref{thm_normality}, like, say,
$h(x)=(1-|x|)_+^\alpha$ for which $\sigma(R,h)\simeq R^{\frac 12-\alpha}$.

\subsection{Few questions}
We do not know whether asymptotic normality holds for all functions
$h\in C^1_0$, or whether the condition
$R^\alpha\sigma(R,h)\to\infty$ is already sufficient for asymptotic
normality of linear statistics associated with a $C^\alpha_0$-function.
Also, we believe that the assertion of Theorem~\ref{thm_normality} can
be extended to functions  $h\in C^{\alpha}\cap L^2_0$ with $-1<\alpha<0$
but our current techniques seem insufficient to handle this case properly.

\subsection{Bounded test-functions}
It is worth mentioning that Theorem~\ref{thm_normality} is complemented
by our recent results~\cite{NS} where, using different methods, we prove that
the fluctuations of $n(R, h)$
are asymptotically normal when $h$ is bounded, and the variance
of $n(R, h)$ grows at least as a positive power of $R$.

\subsection{Test-functions with abnormal fluctuations in linear statistics}
We will see that for every
$\alpha\in(0,1)$, the function $h=|x|^\alpha\psi(x)$, where $\psi$ is a smooth cut-off
that equals $1$ in a neighborhood of the origin, yields an abnormal
behavior of the corresponding linear statistics. Clearly, $h\in C^\alpha_0$ and it
is possible to show that $\sigma(R,h)\simeq R^{-\alpha}$. This
shows that Theorem~\ref{thm_normality} is sharp on a rough power scale.

The reason for the loss of asymptotic normality is that only a small neighbourhood
of the origin containing a bounded number of zeroes of $F$
contributes to the variance of $n(R, h)$, which
is not consistent with the idea of normal fluctuations of linear statistics.

\subsection{Comparing random complex zeroes with limiting eigenvalue process for the Ginibre ensemble}
It is interesting to juxtapose our results with what is known for the $N\to\infty$
limit of the eigenvalue point process $\mathcal G$ of the Ginibre ensemble of $N\times N$
random matrices with independent standard Gaussian complex entries.
The process $\mathcal G$ is a determinantal one,
with the kernel $K(z, w) = e^{z \overline{w}}$ (curiously, this is the same kernel that appears as the
covariance for the G.E.F. $F(z)$). It is known that in this case the variance of linear
statistics never decays, for $W^2_1$-functions $h$ it tends to the limit proportional to $\| \nabla h \|_{L^2}$,
and it never grows faster than $R^2$.
A counterpart of our result proven in ~\cite{NS} is a theorem of Soshnikov. In~\cite{Soshn},
he proved among other things that for arbitrary determinantal point processes, the fluctuations of
linear statistics associated with a compactly supported bounded positive function are normal
if the variance grows at least as a positive power of expectation as the intensity tends to infinity.
A counterpart of the limiting case $\alpha = 2$ in Theorem~\ref{thm_normality} (that is, of the result
from~\cite[Part~I]{ST}) was recently found by Rider and Vir\'ag in \cite{RV}.
They proved that the fluctuations for linear statistics of process $\mathcal G$
are normal
when the test function $h$ belongs to the Sobolev space $W_1^2$. It is not clear whether there is
any meaningful statement interpolating the theorems of Soshnikov and Rider and Vir\'ag. It can happen that
our Theorem~\ref{thm_normality} simply has no interesting counterpart for the process $\mathcal G$.
It is also worth mentioning that the proofs in the determinantal case are quite different from ours.
They are based on peculiar combinatorial identities for the cumulants of linear statistics that are a
special feature of determinantal point processes.

\subsection{Idea of the proof of Theorem~\ref{thm_normality}}
To prove Theorem~\ref{thm_normality}, we follow
a classical idea of S.~N.~Bernstein~\cite{Bern} and approximate the random variable
$n(R, h)$ by a sum of a large number of independent random variables with
negligible error\footnote{The surgery technique developed in Tsirelson's paper~\cite{Tsirel}
is also reminiscent of Bernstein's idea.}. It is worth mentioning that such approximation
becomes possible only after we separate the high and the low frequencies in $h$.
In this approach, independence appears as a
consequence of the almost independence property of G.E.F.'s found in \cite{NSV1, NSV2}.
Roughly speaking, it says that if the compact sets $K_j\subset\R^2$ are well-separated
from each other, then the restrictions of the normalized processes $F^*(z)=F(z)e^{-\frac12|z|^2}$ to
$K_j$  can be simultaneously approximated by the restrictions of normalized independent
realizations of G.E.F.'s with very high probability.

\subsection*{Acknowledgments} We are grateful to Yuri Makarychev who suggested the
idea of the proof of Lemma~\ref{lem_4.2} which is a central part in our proof of
the asymptotic
independence theorem~\ref{thm_almost_indep}, and  to Manjunath Krishnapur
and Boris Tsirelson for very helpful discussions.

\section{The variance}\label{section_variance}

\subsection{Proof of Theorem~\ref{thm_variance}}

We start with several lemmas. First, we find the coefficients of
the Hermite-Ito orthogonal expansion of the function $\log |\zeta|$
in the Hilbert space
$\mathcal H = L^2\left( \C, \tfrac1{\pi} e^{-|\zeta|^2}\right)$.
We use the notation traditional in the probability theory.

Denote by $\mathcal P_m$ the linear subspace of $\mathcal H$ consisting of
algebraic polynomials in $z,\bar z$ of degree at most $m$, and set
$ \mathcal H^{:m:}=\mathcal P_m \ominus P_{m-1} $. By  $:\!\!\zeta^\alpha {\bar\zeta}^\beta\!\!:$
($\alpha, \beta\in\mathbb Z_+$) we denote the orthogonal basis in $\mathcal H$ obtained by
 projecting the polynomials $ \zeta^\alpha \bar\zeta^\beta $ to
$\mathcal H^{:m:}$, $m=\alpha+\beta$. A short computation shows that
$ \|:\!\!\zeta^\alpha {\bar\zeta}^\beta\!\!: \|^2 = \alpha!\beta!$ \cite[Example~3.32]{Janson}.
Since $\log|\zeta|$ is a radial function,
its expansion in this basis contains only the radial polynomials; i.e., the ones with
$\alpha=\beta$. The coefficients of this expansion can be  readily computed:
\begin{lemma}\label{lem_expansion} We have
\[
\log |\zeta| = c_0 + \sum_{\alpha\ge 1} \frac{c_{2\alpha}}{\alpha!} :\!\!|\zeta|^{2\alpha}\!\!:
\]
with $c_0=\Ex\{\log|\zeta|\}=-\tfrac12 \gamma$ ($\gamma$ is the Euler constant) and
$c_{2\alpha}= (-1)^{\alpha+1} \tfrac1{2\alpha}$ for $\alpha\ge 1$.
\end{lemma}

\par\noindent{\em Proof of Lemma~\ref{lem_expansion}:} Suppose $\alpha\ge 1$.
Let $\mathcal H_{\tt rad}$ be the subspace of $\mathcal H$ that consists of
radial functions. If $\Phi\in\mathcal H_{\tt rad}$, then letting $\Phi(\zeta)=\phi(|\zeta|^2)$, we get
\[
\| \Phi \|_{\mathcal H}^2 = \frac1{\pi} \int_0^{2\pi} {\rm d}\theta \,
\int_0^\infty |\phi(r^2)|^2 e^{-r^2}r\,{\rm d}r
= \int_0^\infty |\phi (t)|^2 e^{-t}\, {\rm d}t\,.
\]
This identifies the subspace $\mathcal H_{\tt rad}$ with the space $L^2\left( \R_+, e^{-t} \right)$.
Therefore, $\tfrac1{\alpha!} :\!\! |\zeta|^{2\alpha}\!\!:\, = (-1)^\alpha L_\alpha (|\zeta|^2)$,  where
\[ L_\alpha (t) = \frac1{\alpha!}\, e^t \,
\frac{{\rm d}^\alpha}{{\rm d}t^\alpha}\bigl( t^\alpha e^{-t} \bigr) \]
are Laguerre orthogonal polynomials. In our case,  $\phi (t) = \tfrac12 \log|t|$, whence
\begin{multline*}
c_{2\alpha} = \frac{(-1)^\alpha}{2\alpha!}\,
\int_0^\infty \log t\, \bigl( t^\alpha e^{-t}\bigr)^{(\alpha)}\, {\rm d}t
= -\frac{(-1)^\alpha}{2\alpha!}\, \int_0^\infty \frac1{t}
\bigl( t^\alpha e^{-t}\bigr)^{(\alpha-1)}\, {\rm d}t
\\ = \, ...\, = -\frac{(-1)^\alpha}{2\alpha}\, \int_0^\infty e^{-t}\, {\rm d}t
= \frac{(-1)^{\alpha+1}}{2\alpha}\,,
\end{multline*}
for all $\alpha \ge 1$. For $\alpha=0$, we have
\[
c_0 = \Ex\{\log|\zeta|\} = \frac12\,
\int_0^\infty (\log t) e^{-t}\, {\rm d}t = \frac12 \Gamma'(1) = - \frac12 \gamma\,.
\]
Hence the lemma. \mbox{} \hfill $\Box$

\begin{lemma}\label{lemma_covar}
Suppose $\zeta_i$, $i=1, 2$, are standard complex
Gaussian random variables with $ | \Ex \left\{ \zeta_1 \overline{\zeta_2}\right\} | = \rho $. Then
\[
\Ex \left( \log|\zeta_1| - c_0 \right) \left( \log|\zeta_2| - c_0
\right) = \frac14 \sum_{\alpha \ge 1}
\frac{\rho^{2\alpha}}{\alpha^2}\,,
\]
where $c_0 = \Ex \log |\zeta_i| = -\tfrac12 \gamma$.
\end{lemma}

\par\noindent{\em Proof of Lemma~\ref{lemma_covar}:} Using that
\[
\Ex \left( \frac1{\alpha!} :\!\! |\zeta_1|^{2\alpha} \!\!: \right)
\left( \frac1{\beta!} :\!\! |\zeta_2|^{2\beta}\!\!: \right) =
\begin{cases}
\rho^{2\alpha} &{\rm if}\ \alpha=\beta, \\
0 &{\rm otherwise}
\end{cases}
\]
\cite[Theorem~3.9]{Janson}, we readily derive the result from Lemma~\ref{lem_expansion}. \mbox{}
\hfill $\Box$

\medskip We note that a different
proof of this lemma can be found in~\cite[Lemma~3.3]{SZ}.

\medskip Now, let $f$ be an arbitrary Gaussian analytic function in a domain $G\subset \C$, and let
$\displaystyle f^*(x) = \frac{f(x)}{\sqrt{\Ex |f(x)|^2}}$. Denote by
$\displaystyle n_f = \sum_{a\colon f(a)=0} \delta_a $
the (random) counting measure of its zeroes, and set $\overline{n}_f = n_f - \Ex n_f$. Define the
normalized
$2$-point correlation measure $\nu$ of zeroes of $f$
on $G\times G$ by
$\nu = \Ex \left\{ \overline n_f\times \overline n_f \right\}$. Knowing the two-point measure $\nu$, we easily
recover the variance of the random variable $n(R, h)$:
\[
\V (R, h) = \iint_{\R^2\times \R^2} h\left(\frac{x_1}{R}\right) h\left(\frac{x_2}{R}\right)\,
{\rm d}\nu (x_1, x_2)\,.
\]

\begin{lemma}\label{lemma_2pt_measure}
\[
\nu = \frac1{16\pi^2} \sum_{\alpha\ge 1}
\frac1{\alpha^2} \Delta_{x_1} \Delta_{x_2} \left( \rho^{2\alpha} \right)
\qquad {\rm (as\ a\ distribution)},
\]
where $\rho (x_1, x_2) = \big| \Ex \left( f^*(x_1) \overline{f^*(x_2)}\right) \big|$
is the normalized correlation coefficient of complex Gaussian random variables $f(x_1)$
and $f(x_2)$, and $\Delta_{x_j}$ are (distributional) Laplacians acting on the variables
$x_j$, $j=1, 2$.
\end{lemma}

\par\noindent{\em Proof of Lemma~\ref{lemma_2pt_measure}:} We have
$\overline{n}_f = \tfrac1{2\pi} \Delta \log|f^*|$ (this is the Edelman-Kostlan formula,
see \cite{HKPV, Sodin}), whence
\[
\nu  = \left( \frac1{2\pi} \right)^2 \Delta_{x_1}
\Delta_{x_2} \Ex \left\{ \log |f^*(x_1)| \log |f(x_2^*)| \right\}\,.
\]
Lemma~\ref{lemma_covar} with $\zeta_i = f^*(x_j)$,
$j = 1,2 $, yields (up to terms constant in $x_1$ or
$x_2$ that vanish after we apply the Laplacians in $x_1$ and
$x_2$)
\[
\Ex \left\{ \log |f^*(x_1)| \log |f(x_2^*)| \right\} = \frac14 \sum_{\alpha\ge 1} \frac1{\alpha^2}\,
\rho (x_1, x_2)^{2\alpha} + \big\langle\, {\rm negligible\ terms}\,
\big\rangle,
\]
which implies Lemma~\ref{lemma_2pt_measure}. \mbox{} \hfill $\Box$

\medskip

\par\noindent{\em Proof of Theorem~\ref{thm_variance}:}
It suffices to prove the theorem for $R=1$. The rest readily follows by scaling
\[
x\mapsto R x, \quad h \mapsto h(R^{-1}\,\cdot), \quad \lambda \mapsto R^{-1}\lambda, \quad
\widehat{h} \mapsto R^2\, \widehat{h}(R\,\cdot)\,.
\]

We have
\[
\V (1, h) = \iint_{\R^2\times \R^2} h(x_1) h(x_2)\, {\rm d}\nu (x_1, x_2)\,.
\]
To compute the $2$-point measure $\nu$, we use Lemma~\ref{lemma_covar}. The normalized
correlation coefficient
$\rho (x_1, x_2)$ of the G.E.F. $F$ equals
\[
\rho (x_1, x_2) = \exp\bigl[ -{\rm Re}(x_1\overline{x}_2) - \tfrac12|x_1|^2 - \tfrac12|x_2|^2 \bigr] =
\exp\bigl[ -\tfrac12 |x_1-x_2|^2 \bigr]\,,
\]
whence,
\[
\V (1, h)
= \frac1{16\pi^2} \sum_{\alpha\ge 1} \frac1{\alpha^2}
\iint_{\R^2\times\R^2} h(x_1)h(x_2)\, \Delta_{x_1}\Delta_{x_2} e^{-\alpha |x_1-x_2|^2}\,
{\rm d}A(x_1)\, {\rm d}A(x_2)\,.
\]

Starting from here, we shall assume that $h\in C^2_0$. It doesn't really matter because
it is possible to show that both sides of the last equation are
continuous functionals in $( L^1 \cap L^2 )(\R^2)$.

Now, we use the identity
\[
\iint_{\R^2\times \R^2} K(x_1-x_2) g(x_1) g(x_2)\, {\rm d} A(x_1)\, {\rm d} A(x_2)
= \int_{\R^2} \widehat{K}(\lambda) |\widehat{g}(\lambda)|^2\, {\rm d} A(\lambda)\,,
\]
with $g(x) = \Delta h (x)$ and $K(x)=e^{-\alpha |x|^2}$. With our normalization
of the Fourier transform we have
\[
\widehat{g}(\lambda) = -4\pi^2 |\lambda|^2\widehat{h}(\lambda), \qquad
\widehat{K}(\lambda) = \tfrac{\pi}{\alpha} e^{-\pi^2 |\lambda|^2/\alpha}.
\]
Finally, we get
\[
\V(1, h) =
\int_{\R^2} |\widehat{h}(\lambda)|^2 M( \la )\, {\rm d} A(\la)\,,
\]
with
\[
M(\la) = \pi^3 |\la|^4 \sum_{\alpha\ge 1} \frac1{\alpha^3}\, e^{-\frac{\pi^2}{\alpha} |\la|^2 }
\]
completing the proof of Theorem~\ref{thm_variance}.  \hfill $\Box$

\begin{remark}\label{remark-var}
{\rm Later, we will need the following observation, which can be proven in a similar way.
Denote by $U$ the function $U(x)=\log |F(x)| - \tfrac12 |x|^2$.
Let $g$ be an $L^2(\R^2)$-function
with compact support. Then the variance of the random variable
$ \displaystyle \int_{\R^2} g U {\rm d}A $ is $\lesssim \|g\|_{L^2}^2$.
}
\end{remark}

\subsection{Corollaries}
Later, we will need the following three corollaries to
Theorem~\ref{thm_variance}.

\begin{corollary}\label{cor_Var_A}
For each $h\in (L^1 \cap L^2)(\R^2)$, we have
\[
\V (R, h) \lesssim R^2 \| h\|_{L^2}^2, \qquad  0<R<\infty\,,
\]
and
\[
\V (R, h) = o(R^2), \qquad R\to\infty\,.
\]
\end{corollary}

\medskip
By $W_2^2$ we denote the Sobolev space of $L^2$-functions $h$ with
$\Delta h \in L^2$.

\begin{corollary}\label{cor_Var_B}
For each $h\in L^1\cap W_2^2$,
\[
\V (R, h) \lesssim R^{-2} \| \Delta h \|_{L^2}^2\,, \qquad 0<R<\infty\,.
\]
\end{corollary}

\medskip

\begin{corollary}\label{cor_Var_C}
For any cut-off function $\chi$ with
$|\widehat{\chi} (\lambda)| \lesssim (1+|\la|^2)^{-1}$, we have
\[
\V (R, h) \gtrsim R^{-2} \| \Delta (h * \chi\ci R) \|_{L^2}^2\,, \qquad 0<R<\infty\,.
\]
Here, $\chi\ci R = R^2 \chi (R\,\cdot)$.
\end{corollary}

\section{Almost independence}

The proof of Theorem~\ref{thm_normality} is based on the almost independence
property of G.E.F. introduced in \cite{NSV1, NSV2}. It says that if
$\{ K_j \} $
is a collection of well-separated compact sets in $\C$, then the restrictions
of normalized processes $ F^* \big|_{K_j} $ can be simultaneously
approximated by restrictions $ F_j^* \big|_{K_j} $ of normalized {\em independent} realizations of G.E.F.'s
with very high probability. Here, and everywhere below, $F^*(z)=e^{-\frac12|z|^2} F(z)$.
In \cite{NSV1, NSV2} we used this idea in the case
when $K_j$ were disks. The proofs given in those papers used Taylor expansions of the
shifted G.E.F.'s $F(z+w_j)e^{-z\overline{w}_j}e^{-\frac12 |w_j|^2}$ where $w_j$
were the centers of the disks $K_j$. These proofs cannot
be immediately extended to arbitrary compact sets $K_j$. Here, we give a
more general version of this principle that does not assume anything about
the structure of the sets $K_j$.

\begin{theorem}\label{thm_almost_indep} Let $F$ be a G.E.F.
There exists a numerical constant $A>1$ with the
following property. Given a family of compact sets $K_j$ in $\C$ with diameters
$d(K_j)$, let $\rho_j \ge \sqrt{\log (3+d(K_j))} $. Suppose that
$A\rho_j$-neighbourhoods of the sets $K_j$ are pairwise disjoint. Then
\[
F^* = F_j^* + G_j^* \qquad  on \ K_j,
\]
where $F_j$ are independent G.E.F.'s and
\[
\mathbb P \left\{ \max_{K_j} |G_j^*| \ge e^{-\rho_j^2 }\right \}
\lesssim  \exp\bigl[ - e^{\rho_j^2} \bigr] \,.
\]
\end{theorem}

When proving this theorem, it will be convenient to treat complex Gaussian random variables
as elements in a ``big'' Gaussian Hilbert space $\mathcal H$. First, we consider the
Gaussian Hilbert space $\mathcal H_F$, which is a closure of finite linear combinations
$U = \sum_k c_k F(z_k) $ with respect to the scalar product generated by the covariance:
$\langle U, V \rangle = \Ex \bigl\{ U \overline{V} \bigr\} $. We assume that the big space
$\mathcal H$ consists of complex valued Gaussian random variables and
contains countably many mutually orthogonal copies of $\mathcal H_F$. This will
allow us to introduce new independent copies of some Gaussian random variables
when necessary.

The proof of Theorem~\ref{thm_almost_indep} goes as follows. First, for each compact set $K_j$,
we choose a sufficiently dense net $Z_j$ and consider the bunch $N_j = \bigl\{ v_z\colon z\in Z_j \bigr\}$
of unit vectors $v_z = F^*(z)$. Since the compact sets
$K_j$ are well-separated, the bunches $N_j$ are almost orthogonal to each other. Then
we slightly perturb the vectors $v_z$ without changing the angles between
the vectors within each bunch $N_j$, making the bunches orthogonal to each other.
More accurately, we construct new bunches
$\widetilde{N}_j = \bigl\{ \widetilde{v}_z\colon z\in Z_j \bigr\}$ so that for $z\in Z_j$,
$\zeta\in Z_k$,
\[
\langle \widetilde{v}_z, \widetilde{v}_\zeta\rangle
=
\begin{cases}
\langle v_z, v_\zeta\rangle & {\rm for\ } j=k, \\
0 & {\rm for\ } j\ne k
\end{cases}
\]
with good control of the  errors
$ \| v_z - \widetilde{v}_z\| $.
Then we extend the Gaussian bunches
$\bigl\{ \widetilde{v}_z e^{\frac12 |z|^2}\colon z\in Z_j \bigr\}$
to {\em independent} G.E.F.'s $F_j$. The difference $G_j = F-F_j$ is a
random entire function
that is small on the net $Z_j$ with probability very close to one.
At the last step of the proof,
we show that $G_j^*$ is small everywhere on $K_j$.

\medskip\par\noindent{\em Proof of Theorem~\ref{thm_almost_indep}:} We start with
the construction of the net $Z_j$. We fix a unit lattice in $\R^2$
and denote by $\mathcal C_j$ a finite collection of the lattice points that
is contained in a $ \tfrac1{\sqrt{2}} $-neighbourhood of the compact set $K_j$ and
is a $\tfrac1{\sqrt 2}$-net for $K_j$. Then we
consider the collection $\mathcal D_j$ of the unit disks centered at $\mathcal C_j$ and
choose about $A^2 \rho_j^2$ equidistant points on the boundary of each of these unit disks
where $A$ is a sufficiently big constant to be chosen later.
The collection of all these points
will be our net $Z_j$. Then,  for $z\in Z_j$, we introduce the bunch
$N_j = \bigl\{ v_z\colon z\in Z_j \bigr\} \subset \mathcal H$ of unit vectors $v_z = F^*(z)$.
Since
\begin{equation}\label{eq_4.1}
\left| \langle v_z, v_\zeta \rangle \right| \le e^{-|z-\zeta|^2/2}
\end{equation}
and since the compact sets $K_j$ are well separated,
the vectors in the bunches $N_j$ and $N_k$ are almost orthogonal when $j\ne k$.

\smallskip Now, we replace almost orthogonal bunches by the orthogonal ones.
\begin{lemma}\label{lem_4.2} There exist vectors $w_z\in\mathcal H$ such that
\begin{equation}\label{eq_4.5}
\langle w_z, w_\zeta \rangle =
\begin{cases}
- \langle v_z, v_\zeta \rangle, & z\in Z_j, \ \zeta\in Z_k, \ j\ne k, \\
0, & z, \zeta \in Z_j, \ z\ne \zeta, \\
e^{-\frac15 A^2 \rho_j^2}, & z=\zeta\,,
\end{cases}
\end{equation}
and
\begin{equation}\label{eq_4.6}
\langle w_z, v_\zeta \rangle = 0 \qquad {\rm for\ each \ } z,\zeta\,.
\end{equation}
\end{lemma}

\par\noindent{\em Proof of Lemma~\ref{lem_4.2}:} Consider the Hermitian matrix $\Gamma$ with
the elements $\gamma (z, \zeta)$ defined by the right-hand side of \eqref{eq_4.5}. We will check that
the matrix $\Gamma$ is positive definite. This will ensure the existence of vectors $w_z$ in
$\mathcal H \ominus {\rm span}(v_z)$ with the Gram matrix $\Gamma$.

By the classical Gershgorin theorem, each eigenvalue of the matrix $\Gamma$ lies in
one of the intervals
\[
\bigl( \gamma (z,z) - t(z), \gamma (z, z) + t(z) \bigr) \qquad{\rm with} \quad
t(z) = \sum_{\zeta\colon \zeta\ne z} |\gamma (z, \zeta)|.
\]
Below, we check that for each $z$,
$t(z) < \gamma (z,z)$; i.e, all eigenvalues of the matrix $\Gamma$ are positive, whence
$\Gamma > 0$.

\medskip
\begin{lemma}\label{lem_4.1}
Let $z\in Z_j$. Then
\begin{equation}\label{eq_4.4}
\sum_{{k\colon k\ne j}}\, \sum_{\zeta\in Z_k}
\left|  \langle v_z, v_\zeta \rangle \right| < e^{- \frac15 A^2 \rho_j^2}\,,
\end{equation}
provided that the constant $A$ in the assumptions of Theorem~\ref{thm_almost_indep}
is big enough.
\end{lemma}

\par\noindent{\em Proof of Lemma~\ref{lem_4.1}:} Fix $k\ne j$ and $c\in\mathcal C_k$.
Denote by $D$ the unit disk centered at $c$. Then for each point
$\zeta\in\partial D$, we have
$|z-\zeta|\ge |z-c| -1$, and by \eqref{eq_4.1}
\[
\left| \langle v_z, v_\zeta \rangle \right| \le e^{-\frac14 |z-c|^2 + 1}\,,
\]
whence
\[
\sum_{\zeta\in \partial D_w }\left| \langle v_z, v_\zeta \rangle \right| \le
A^2 \rho_k^2\, e \cdot e^{-\frac14 |z-c|^2}\,,
\]
and
\[
\sum_{\zeta\in Z_k }\left| \langle v_z, v_\zeta \rangle \right| \le
A^2 e \sum_{c\in \mathcal C_k} \rho_k^2 e^{-\frac14 |z-c|^2}\,.
\]
Since the $A\rho_j$-neighbourhood of $K_j$ and the $A\rho_k$-neighbourhood of
$K_k$ are disjoint, for each $z\in Z_j$ and $c\in\mathcal C_k$
with $j\ne k$, we have
\[
\rho_k < \rho_j + \rho_k \le \frac1{A} \bigl( |z-c| + \text{dist}(z, K_j) + \text{dist}(c, K_k) \bigr)
\le \frac1{A} \left( |z-c| + \sqrt 2 + 1 \right).
\]
Since $\rho_k \ge 1$ and $A$ is sufficiently big, we conclude that
$ \rho_k < \tfrac2{A}|z-c| $. Then
\[
\sum_{\zeta\in Z_k }\left| \langle v_z, v_\zeta \rangle \right| \le
4e \sum_{c\in \mathcal C_k} |z-c|^2 e^{-\frac14 |z-c|^2}\,,
\]
and
\[
\sum_{{k\colon k\ne j}} \, \sum_{\zeta\in Z_k }\left| \langle v_z, v_\zeta \rangle \right| \le
4e \sum_{k\ne j }\sum_{c\in \mathcal C_k} |z-c|^2 e^{-\frac14 |z-c|^2}\,.
\]

To estimate the right-hand side, we introduce the counting measure $\mu$ of the
set $\displaystyle \bigcup_{{k\colon k\ne j}} \mathcal C_k$. By our construction, for $z\in Z_j$, we have
\[
\mu (D(z, t)) \le
\begin{cases}
C_1 t^2, & \text{for\ } t\ge A\rho_j - C_2, \\
0, & \text{for\ } t< A\rho_j - C_2,
\end{cases}
\]
with positive numerical constants $C_1$ and $C_2$.
Then
\begin{multline*}
4e \sum_{k\ne j }\sum_{c\in \mathcal C_k} |z-c|^2 e^{-\frac14 |z-c|^2/2}
= 4e \int_{\R^2} |z-c|^2 e^{-\frac14 |z-c|^2}\, {\rm d}\mu (c) \\
= 4e \int_0^\infty (-2t+\tfrac12 t^3) e^{-t^2/4} \mu (D(z, t)) \, {\rm d} t \\
\le C_3 \int_{A\rho_j-C_2}^\infty t^5 e^{-t^2/4}\, {\rm d} t
<  e^{- \frac15 A^2 \rho_j^2}\,,
\end{multline*}
provided that the constant $A$ is big enough. This proves both lemmas. \hfill $\Box$

\medskip We resume the proof of Theorem~\ref{thm_almost_indep} and
denote by $\mathcal H_j$ the linear span in $\mathcal H$ of the vectors
$\left\{ v_z+w_z\colon z\in Z_j \right\}$ where $w_z$ are the vectors from
Lemma~\ref{lem_4.2}. By construction, the subspaces $\mathcal H_j$ are mutually orthogonal.
Let $P_{\mathcal H_j}$ be the orthogonal projection to $\mathcal H_j$ and let
$u_z = v_z - P_{\mathcal H_j} v_z$. Then
\[
\| u_z \| \le \| w_z \| \le e^{-\frac1{10} A^2 \rho_j^2}\,.
\]
Now, we choose subspaces $\mathcal H_j^*$ orthogonal to all subspaces $\mathcal H_j$ and orthogonal
to each other, and choose the isometries $S_j$ that move all the vectors $u_z$ with $z\in Z_j$
to $\mathcal H_j^*$. At last, we set
\[
\widetilde{v}_z = P_{\mathcal H_j} v_z + S_j u_z\,,
\]
and introduce the new bunches of unit vectors
$\widetilde{N}_j = \left\{ \widetilde{v}_z\colon z\in Z_j \right\}$. The new bunches are mutually orthogonal.
Indeed, if $z\in Z_j$ and $\zeta\in Z_k$ with $k\ne j$, then
\[
\langle \widetilde{v}_z, \widetilde{v}_\zeta \rangle
= \langle P_{\mathcal H_j} v_z + S_j u_z , P_{\mathcal H_k} v_\zeta + S_k u_\zeta \rangle = 0
\]
since the subspaces $\mathcal H_j$, $\mathcal H_j^*$, $\mathcal H_k$, $\mathcal H_k^*$ are mutually
orthogonal for $k\ne j$. On the other hand, the angles between the vectors within the same bunch remain
the same: for $z, \zeta \in Z_j$, we have
\[
\langle \widetilde{v}_z, \widetilde{v}_\zeta \rangle
= \langle P_{\mathcal H_j} v_z + S_j u_z , P_{\mathcal H_j} v_\zeta + S_j u_\zeta \rangle
= \langle P_{\mathcal H_j}v_z , P_{\mathcal H_j} v_\zeta \rangle
+ \langle u_z , u_\zeta \rangle = \langle v_z, v_\zeta \rangle\,.
\]
Hence the vectors $\left\{ \widetilde{v}_z \right\}$ define the same Gaussian process on $Z_j$ as
$\left\{ v_z \right\}$, and the new processes are independent.

We set $M_j = e^{\frac1{20} A^2 \rho_j^2}$ and note that
$ \| \widetilde{v}_z - v_z \| \le 2 \| u_z \| \le 2M_j^{-2} $
for all $z\in Z_j$.  Then
\[
\mathbb P \bigl\{  |\widetilde{v}_z - v_z| > t  \bigr\} = e^{ -(t/\| \widetilde{v}_z - v_z \|)^2 }
\le e^{-(tM_j^2/2)^2}\,, \qquad z\in Z_j\,.
\]
Plugging in this estimate with $ t=\tfrac1{M_j} $, we get
\[
\mathbb P \Big\{  \sup_{z\in Z_j } |\widetilde{v}_z - v_z| > \tfrac1{M_j} \Big\}
\le \# Z_j \cdot e^{-\frac14 M_j^2}\,.
\]
Since
\[
\# Z_j \le A^2 \rho_j^2 \cdot \# \mathcal C_j \lesssim
A^2 \rho_j^2 \cdot (1+d(K_j))^2  \lesssim A^2 \rho_j^2 \cdot e^{2\rho_j^2}\,,
\]
we have
\[
\mathbb P \Big\{  \sup_{z\in Z_j } |\widetilde{v}_z - v_z| > \tfrac1{M_j} \Big\}
\lesssim  A^2 \rho_j^2 \cdot e^{2\rho_j^2} \cdot e^{-\frac14 M_j^2} < e^{-cM_j^2}\,,
\]
provided that $A$ is sufficiently big.

Now, we extend the processes $\left\{ \widetilde{v}_z e^{\frac12 |z|^2}\colon z\in Z_j \right\}$
to {\em independent} G.E.F.'s $F_j$. We know that $G_j^* = F^* - F_j^*$ is small on the net $Z_j$
with probability very close to $1$:
\begin{equation}\label{eq-missed}
\mathbb P \Big\{  \sup_{z\in Z_j} |G_j^*(z)| > \tfrac1{M_j} \Big\}
< e^{-cM_j^2}\,.
\end{equation}
Discarding the event $ \bigl\{  \sup_{z\in Z_j} |G_j^*(z)| > M_j^{-1} \bigr\} $, we assume that
\begin{equation}\label{eq-assumption}
\sup_{z\in Z_j} |G_j^*(z)| \le \frac1{M_j}
\end{equation}
and estimate the difference $G_j^* = F^* - F_j^*$ everywhere on $K_j$.

\medskip Let $\mathcal D $ be one of the unit disks with the center $\kappa\in\mathcal C_j$. It
is convenient to move the point $\kappa$ to the origin using the projective translation
\[
(T_\kappa f)(z) = f(\kappa +z) e^{-z\overline{\kappa}} e^{-\frac12|\kappa|^2}\,.
\]
Comparing the covariances, it is easy to check that if $F$ is a G.E.F., then
$T_\kappa F$ is also a G.E.F.~\cite[Lemma~2.6]{NSV1}. To simplify our notation, we set
$ G(z) = \bigl( T_\kappa G_j \bigr)(z) $,
and denote by $ \Lambda $ the set of equidistant points on the unit circumference such that
$\Lambda+\kappa \subset Z_j$. Due to our assumption~\eqref{eq-assumption},
the random analytic function $G$ is very small on $\Lambda$:
\begin{equation}\label{eq-estG1}
|G(\la)| = | T_\kappa (F-F_j)(\la) | \le
\tfrac1{M_j} e^{\frac12|\la+\kappa|^2 - \frac12 |\kappa|^2 - {\rm Re}(\la \overline{\kappa})}
= \tfrac1{M_j} e^{\frac12} < \tfrac2{M_j}
\end{equation}
Since $ |G|\le |T_\kappa F| + |T_\kappa F_j| $,
the function $G$ is not too big everywhere on the disk $\{ |z|\le 2 \} $  with very high probability:
\begin{equation}\label{eq-estG2}
\mathbb P \left\{ \max_{\{|z|\le 2\}} |G (z)| \ge M_j \right\} \le e^{-cM_j^2}\,.
\end{equation}
This follows from Lemma~6 in \cite{NSV2}
applied to the G.E.F.'s $T_\kappa F$ and $T_\kappa F_j$.

Now, we suppose that $ |G|\le M_j $ everywhere in the disk $\{|z|\le 2\}$.
Then estimate~\eqref{eq-estG1} yields that $G$ is very small everywhere on the
disk $\{ |z| \le \tfrac1{\sqrt 2}\}$. To see this, first, we approximate $G$ by its Taylor
polynomial $P_{N}$ of sufficiently large degree $N$. Let
\[
G(z) = \sum_{n\ge 0} g_{n} z^n\,.
\]
By Cauchy's inequalities, for each $n\in\mathbb Z_+$, $|g_n| \le M_j\cdot 2^{-n}$,
and
\[
\left| \sum_{n\ge N+1} g_{n} z^n \right| \le \frac{M_j}{2^N}\,.
\]
Since  we wish the approximation error $\max_D |G - P_{N}|$ to be less
than $M_j^{-1}$, we need to take $2^N \ge M_j^2$, that is,
$ N \approx 2 \log_2 M_j $. Then $|P_{N}| \le 2M_j^{-1}$ at about
$A^2 \rho_j^2 = 20 \log M_j$
equidistant points $\la$ on the unit circumference. Since the number of these points is
bigger than the degree of $P_{N}$, and since these points are equidistant, we have
\[
g_n = \frac1{\# \Lambda} \sum_{\la\in\Lambda} P_N(\la) \la^{-n}, \qquad 0 \le n \le N,
\]
whence, $|g_n|\le 2M_j^{-1}$ for $0\le n \le N$. Thus,
\[ |P_N (z)| < \frac2{M_j}\, \frac{\sqrt 2}{\sqrt{2}-1} < \frac8{M_j}, \]
and $|G(z)|< \tfrac9{M_j}$ for $|z|\le \tfrac1{\sqrt 2}$.

It is easy to check that $ |G_j^*(z)| = |G(z-\kappa)| e^{-\frac12|z-\kappa|^2} $. Hence,
for $ z\in D (\kappa, \tfrac1{\sqrt 2}) $, we have $ |G_j^*(z)| < \tfrac9{M} $,
provided that $|G(z)|< \tfrac9{M_j}$ therein.
Since the union of all disks of radius $\tfrac1{\sqrt 2}$ centered at $\mathcal C_j$
covers $K_j$, we see that
\begin{multline*}
\mathbb P \bigl\{ \max_{K_j} |G_j^*| > e^{-\rho_j^2} \bigr\}
< \mathbb P \bigl\{ \max_{K_j} |G_j^*| > 9M_j^{-1} \bigr\} \\
\stackrel{\eqref{eq-missed}, \eqref{eq-estG2}}< C \cdot \# \mathcal C_j \cdot e^{-cM_j^2}
= C \exp\bigl[ 2\rho_j^2 - ce^{\frac1{10} A^2 \rho_j^2} \bigr] < \exp \bigr[ -e^{\rho_j^2} \bigl],
\end{multline*}
provided that the constant $A$ is big enough. This
completes the proof of the almost independence theorem. \hfill $\Box$

\section{Asymptotic normality. Proof of Theorem~\ref{thm_normality}}\label{section_AsymptNormality_I}

Without loss of generality, we assume till the end of the proof that $\alpha<2$.
We fix a $C^\alpha_0$-function $h$ supported by the square $[-\tfrac12, \tfrac12]\times [-\tfrac12, \tfrac12]$
and fix a big positive $ R $ such that $\sigma (R, h) \gtrsim R^{-\alpha+\epsilon}$.
By $ \mathcal S $ we denote the square $ [-1, 1] \times [-1, 1] $.
As above, we put $\overline{n}(R, h) = n(R, h) - \Ex n(R, h)$.

\subsection{Preliminary smoothing} First, we approximate the function
$h$ by a $C^2_0$-function $h_R$ so that
\begin{equation}\label{eq_prelim_1}
\| h-h_R \|_\infty \lesssim R^{-3}\,,
\end{equation}
\begin{equation}\label{eq_prelim_2}
\| \nabla h_R \|_\infty\,, \| \Delta h_R\|_\infty \lesssim R^{M(\alpha)}\,,
\end{equation}
\begin{equation}\label{eq_prelim_3}
\| h_R \|_{C^\alpha} \lesssim \| h \|_{C^\alpha}\,.
\end{equation}
For instance, we may take $h_R = h * \phi_\epsilon * \phi_\epsilon$, where
$\phi_\epsilon = \tfrac1{\pi \epsilon^2} \done_{\{ |x|<\epsilon \}}$ with
$\epsilon = R^{-3/\alpha}$. This gives us \eqref{eq_prelim_2} with
$M(\alpha) = \tfrac6{\alpha} $.

Now, $n(R, h) = n(R, h_R) + n(R, h-h_R)$, and
\[
\V (R, h-h_R) \stackrel{\rm Cor~\ref{cor_Var_A}}\lesssim R^2 \| h-h_R\|_{L^2}^2 \lesssim R^{-4}
= o(1) \V(R, h)
\]
since by \eqref{eq_variance}, $V(R, h)$ cannot decay faster than $R^{-2}$.
Hence, when proving asymptotic normality of the linear statistics $n(R, h)$, we may replace it
by $n(R, h_R)$. To simplify the notation, we omit the subscript and continue to denote the
function $h_R$ by $h$.

\subsection{Separating low and high frequencies}

We fix a radial function $\chi\in C^\infty_0$ with
\begin{equation}\label{eq_FTchi}
\left| \widehat{\chi} (\la) - 1 \right| =O(|\la|^2)\,, \qquad \la \to 0
\end{equation}
and decompose the function $h$ into low and high frequency
parts as follows
\[
h = h*\chi\ci R + (h-h*\chi\ci R) = h_{\tt L} + h_{\tt H}\,.
\]
Here, as above, $\chi\ci R = R^2 \chi (R\, \cdot \,)$. We put $g = \Delta h_{\tt L}$.

\begin{lemma}\label{lemma_estimate_h_g}
We have
\begin{equation}\label{eq_h}
\| h_{\tt H} \|_\infty \lesssim \|h\|_{C^\alpha} R^{-\alpha}\,,
\end{equation}
and
\begin{equation}\label{eq_g}
\| g \|_\infty \lesssim \|h\|_{C^\alpha} R^{2-\alpha}\,.
\end{equation}
\end{lemma}

\par\noindent{\em Proof of Lemma~\ref{lemma_estimate_h_g}:} We give the proof only in the
case $\alpha>1$, the proof in the other case $\alpha\le 1$ is very similar.
Since the function $\chi$ is radial, we
have
\[
h_{\tt H}(x) = \int_{\R^2} h(x+y)  \chi\ci R (y)\, {\rm d} A(y) \\
= \int_{\R^2} \left[ h(x+y) - h(x) - \langle \nabla h(x), y \rangle \right] \chi\ci R (y)\, {\rm d} A(y)
\]
(we've used that $\chi\ci R$ is orthogonal to linear functions).
Therefore,
\[
\| h_{\tt H} \|_\infty \lesssim \| h \|_{C^\alpha} \int_{\R^2} |y|^\alpha |\chi\ci R (y)| \, {\rm d} A(y)
\lesssim \|h\|_{C^\alpha} R^{-\alpha}\,.
\]
Next,
\[
g(x) = \int_{\R^2} h(x+y) \Delta \chi\ci R(y)\, {\rm d} A(y)
= \int_{\R^2} \left[ h(x+y) - h(x) - \langle \nabla h(x), y \rangle
\right] \Delta \chi\ci R(y)\, {\rm d} A(y)\,,
\]
whence
\[
\| g \|_\infty \lesssim \|h\|_{C^\alpha}  \int_{\R^2} |y|^\alpha
| \Delta \chi\ci R(y) |\, {\rm d} A(y) \lesssim \|h\|_{C^\alpha} R^{2-\alpha}\,,
\]
completing the proof. \hfill $\Box$

\begin{lemma}\label{lemma_L2estimate_h_g}
We have
\begin{equation}\label{eq_L2_h}
\| h_{\tt H} \|_{L^2}^2 \lesssim R^{-2} \V (R, h) \,,
\end{equation}
and
\begin{equation}\label{eq_L2_g}
\| g \|_{L^2}^2 \lesssim R^2 \V (R, h)\,.
\end{equation}
\end{lemma}

\par\noindent{\em Proof of Lemma~\ref{lemma_L2estimate_h_g}}:
We have
\[
R^2 \|h_{\tt H}\|_{L^2}^2 = R^2 \int_{\R^2} \bigl| \widehat{h}(\la) \bigr|^2 \cdot
\bigl| 1 - \widehat{\chi}(\tfrac1R \la )\bigr|^2\, {\rm d} A(\la)\,.
\]
Since the Fourier transform of the cut-off function $\chi$ satisfies \eqref{eq_FTchi},
Theorem~\ref{thm_variance} yields \eqref{eq_L2_h}. The second estimate
\eqref{eq_L2_g} follows from Corollary~\ref{cor_Var_C}. \hfill $\Box$

\medskip Now, we split the linear statistics ${\overline n}(R, h)$ into low and high frequency parts:
\begin{multline*}
\overline{n}(R, h) = \frac1{2\pi R^2}
\int_{R\, \mathcal S} g\left( \frac{x}{R} \right) U(x)\, {\rm d} A(x)
+ \frac1{2\pi}
\int_{R\, \mathcal S} h_{\tt H} \left( \frac{x}{R} \right) \Delta U(x)\, {\rm d} A(x) \\
= \overline{n}_{\tt L}(R, h) + \overline{n}_{\tt H}(R, h)\,.
\end{multline*}
Here, as before, $U = \log |F^*|$ is the random potential.
Note that
$\Ex \overline{n}_L = 0$ (since $g=\Delta h_{\tt L}$, the integral of $g$ against constants
vanishes), and therefore $\Ex \overline{n}_{\tt H} = 0$.

\subsection{Squares and corridors}

We fix the parameters $0<\beta<\gamma < \tfrac12 \epsilon $ with $\epsilon$ taken from
the assumptions of the theorem,
and partition the plane into the squares $\mathcal Q$ with sides
parallel to the coordinate axes and of length $R^{\gamma-1}$ separated by the corridors
of width $R^{\beta-1}$. We denote by $K$ the union of the corridors.
The partition is defined up to translations of the position of the system
of corridors. We use this freedom to discard the contribution of the corridors.
Averaging over translations,
it is not difficult to see that there exists a position of the corridors such that
\begin{equation}\label{eq_g1}
\| g \done_{K}\|_{L^2}^2 \lesssim R^{-(\gamma-\beta)} \|g\|_{L^2}^2\,,
\end{equation}
and
\begin{equation}\label{eq_h1}
\| h_{\tt H} \done_{K}\|_{L^2}^2
\lesssim R^{-(\gamma-\beta)} \|h_{\tt H}\|_{L^2}^2\,.
\end{equation}
Then we immediately obtain
\begin{lemma}\label{lem_channels}
There exists a position of the corridors such that
\begin{multline*}
\V \left\{ \frac1{R^2} \int_{R K} g\left( \frac{x}{R} \right) U(x)\, {\rm d} A(x)  \right\}
+ \V \left\{ \int_{R K} h_{\tt H}\left( \frac{x}{R} \right) \Delta U(x)\, {\rm d} A(x)  \right\}
\\
= o(1) \V (R, h)\,, \qquad R\to\infty\,.
\end{multline*}
\end{lemma}

\medskip\par\noindent{\em Proof of Lemma~\ref{lem_channels}:}
We fix position of the corridors to satisfy \eqref{eq_g1} and \eqref{eq_h1}.
It follows from Remark~\ref{remark-var}
that for any test function $g$, the variance of the random variable
$\int_{\R^2} g(x) U(x)\, {\rm d} A(x)$ is $\lesssim \| g \|_{L^2}^2$.
Scaling by $R$, we see that
\begin{multline*}
\V \left\{ \frac1{R^2} \int_{R K} g\left( \frac{x}{R} \right) U(x)\, {\rm d} A(x)  \right\}
\lesssim \frac1{R^4} \cdot R^2 \| g  \done_{K}\|_{L^2}^2 \\
\stackrel{\rm\eqref{eq_g1}}\lesssim R^{-2-(\gamma-\beta)} \|g\|_{L^2}^2
\stackrel{\rm\eqref{eq_L2_g}}= o(1) \V (R, h)\,.
\end{multline*}
Similarly,
\begin{multline*}
\V \left\{ \int_{R K} h_{\tt H}\left( \frac{x}{R} \right) \Delta U(x)\, {\rm d} A(x) \right\}
\stackrel{\rm Cor~\ref{cor_Var_B}}\lesssim R^2 \| h_{\tt H} \done_{K} \|_{L^2}^2 \\
\stackrel{\rm\eqref{eq_h1}}\lesssim R^{-(\gamma-\beta)}\, R^2 \|h_{\tt H}\|_{L^2}^2
\stackrel{\rm\eqref{eq_L2_h}}= o(1) \V (R, h) \,.
\end{multline*}
This proves the lemma. \hfill $\Box$

\subsection{Smoothing the indicator-functions of the squares $\mathcal Q$}
Later, when we will deal with the random variables
\[
\int_{R \mathcal Q} h_{\tt H}\left( \frac{x}{R}\right)
\Delta U(x)\, {\rm d} A(x)
\]
it will be convenient to modify the functions $h_{\tt H} \done_{\mathcal Q}$
smoothing them near the boundaries $\partial\mathcal Q$ with the ``range of smoothing''
about $R^{\beta-1}$. We introduce a `smooth indicator function' $\widetilde{\done}_{\mathcal Q}$
so that $ 0 \le \widetilde{\done}_{\mathcal Q} \le 1$,
\[
\widetilde{\done}_{\mathcal Q} =
\begin{cases}
1 & {\rm on\ } \mathcal Q\,, \\
0 & {\rm when\ } {\rm dist}(x, \mathcal Q) \ge \tfrac1{10} R^{\beta-1}\,,
\end{cases}
\]
$\| \nabla \widetilde{\done}_{\mathcal Q} \|_\infty \lesssim R^{1-\beta}$,
and $\| \Delta \widetilde{\done}_{\mathcal Q} \|_\infty \lesssim R^{2(1-\beta)}$.
Let $\widetilde{\mathcal Q}
= {\rm supp}(\widetilde{\done}_{\mathcal Q})$,
$h_{\tt H}^{\mathcal Q} = h_{\tt H}\widetilde{\done}_{\mathcal Q} $.

\begin{lemma}\label{lem_FuzzySq} We have

\smallskip\par\noindent {\rm (i)} $ h_{\tt H}^{\mathcal Q}(x) = h_{\tt H}(x) $ on  $\mathcal Q$;

\smallskip\par\noindent {\rm (ii)} $\| h_{\tt H}^{\mathcal Q} \|_\infty
\lesssim \| h\|_{C^\alpha} R^{-\alpha}$;

\smallskip\par\noindent {\rm (iii)}
$ \| \Delta h_{\tt H}^\mathcal Q \|_\infty \lesssim C(h)\, R^{M(\alpha)+2 } $;

\smallskip\par\noindent {\rm (iv)}
\[
\V \left\{ \sum_{\mathcal Q } \int_{\R^2}
\left( h_{\tt H}\done_{\mathcal Q} - h_{\tt H}^{\mathcal Q}\right)\left( \frac{x}{R} \right)
\Delta U(x)\, {\rm d} A(x) \right\} \lesssim R^{-(\gamma - \beta)} \V (R, h)\,.
\]
\end{lemma}

\par\noindent{\em Proof of Lemma~\ref{lem_FuzzySq}:} The properties (i)-(iii) are
straightforward. Property (iv) follows from
\begin{multline*}
\V \left\{ \sum_{\mathcal Q } \int_{\R^2}
\left( h_{\tt H} \done_{\mathcal Q} - h_{\tt H}^{\mathcal Q}\right)\left( \frac{x}{R} \right)
\Delta U(x)\, {\rm d} A(x) \right\} \\
\stackrel{\rm Cor~\ref{cor_Var_B}}\le
R^2 \Bigl\|\, \sum_{\mathcal Q} \bigl( h_{\tt H} \done_{\mathcal Q} - h_{\tt H}^{\mathcal Q}
\bigr)\, \Bigr\|^2_{L^2}
\lesssim R^2 \| h_{\tt H} \done_{K}\|^2_{L^2} \\
\stackrel{\rm \eqref{eq_h}}\lesssim R^{-(\gamma-\beta)}\, R^2 \| h_{\tt H}\|^2_{L^2}
\stackrel{\eqref{eq_L2_h}}\lesssim R^{-(\gamma - \beta)} \V (R, h)\,.
\end{multline*}
This proves the lemma. \hfill $\Box$

\medskip

Now, we introduce the random variables
\[
\xi^{\mathcal Q}_{\tt L}(R) = \frac1{2\pi R^2} \int_{R \mathcal Q} g \left( \frac{x}{R} \right)
U(x)\, {\rm d} A(x)\,,
\]
and
\[
\xi^{\mathcal Q}_{\tt H}(R) = \frac1{2\pi } \int_{R \mathcal Q} h_{\tt H}^{\mathcal Q}\left( \frac{x}{R} \right)
\Delta U(x)\, {\rm d} A(x)\,,
\]
and set
\begin{multline*}
\Theta (R) \stackrel{\rm def}=
\frac1{2\pi R^2} \int_{R K} g\left( \frac{x}{R} \right) U(x)\, {\rm d} A(x)
+ \frac1{2\pi} \int_{R K} h_{\tt H}\left( \frac{x}{R} \right) \Delta U(x)\, {\rm d} A(x) \\
+ \frac1{2 \pi} \sum_{\mathcal Q } \int_{\R^2}
\left( h_{\tt H} \done_{\mathcal Q} - h_{\tt H}^{\mathcal Q}\right)\left( \frac{x}{R} \right)
\Delta U(x)\, {\rm d} A(x)\,.
\end{multline*}
Then
\[
\overline{n}(R, h) = \Theta (R) + \sum_{\mathcal Q }  \underbrace{ \left( \xi^{\mathcal Q}_{\tt L} (R) +
\xi^{\mathcal Q}_{\tt H} (R) \right)}_{\stackrel{\rm def}=  \xi^{\mathcal Q} (R)} \,.
\]

\subsection{Probabilistic lemma}
To prove the asymptotic normality of the random variables $\overline{n}(R, h)$ we use the
following lemma:

\begin{lemma}\label{lemma_prob} Fix $p\ge 4$.
Assume that $\xi$ is a real random variable with variance $\sigma^2$ such that
\[
\xi = \Theta + \sum_{j=1}^N \xi_j
\]
where

\smallskip\par\noindent 1) $\V \left\{ \Theta \right\} \le \delta \sigma^2$;

\smallskip\par\noindent 2) $\sum_j \Ex |\xi_j|^p \le \delta \sigma^p$;

\smallskip\par\noindent 3) there exist independent random variables $\eta_j$ such that
\[
\sum_j \mathbb P \left\{ |\xi_j - \eta_j| > \delta N^{-1} \sigma \right\} \le \delta N^{-2}\,.
\]
Then the random variable $\sigma^{-1} \left( \xi - \Ex \xi\right)$ is close in distribution to the normal law
when $\delta$ is close to $0$.
\end{lemma}

Due to the second condition, if $\delta\to 0$, then automatically
$N\to\infty$.

\medskip\par\noindent{\em Proof of Lemma~\ref{lemma_prob}:} using appropriate truncations,
we will reduce this lemma to one of
the standard versions of the central limit theorem.

Without loss of generality, we assume that $\sigma=1$. Replacing
$\delta$ by $2^p \delta$, we also assume that $\Ex \xi =
\Ex \Theta = \Ex \xi_j = 0$.
Put $\eta_j' = \sign \eta_j \cdot \min\left\{ |\eta_j|, N^{2/p} \right\}$ and denote
\[
\Omega_j = \left\{ |\xi_j - \eta_j| > \delta N^{-1} \right\}, \qquad
\Omega_j' = \left\{ |\xi_j - \eta_j'| > \delta N^{-1} \right\}\,.
\]
We have
\[
\mathbb P \left( \Omega_j' \right) \le \mathbb P (\Omega_j) + \mathbb P \left( \{|\xi_j|>N^{2/p}\} \right)
\le \mathbb P (\Omega_j) + N^{-2} \Ex |\xi_j|^p\,,
\]
whence
\[
\sum_j \mathbb P \left( \Omega_j' \right) \le \delta N^{-2} + N^{-2} \delta = 2\delta N^{-2}\,.
\]
Next,
\begin{multline*}
\Ex |\xi_j - \eta_j'|^2 \le \delta^2 N^{-2} + 2 \int_{\Omega_j'} |\xi_j|^2\, {\rm d} \mathbb P
+ 2 \int_{\Omega_j'} |\eta_j'|^2\, {\rm d} \mathbb P \\
\le \delta^2 N^{-2} + 2 \int_{\Omega_j'} \left[ N^{-\frac2{p} (p-2)} |\xi_j|^p + N^{\frac4{p}}
\right] \, {\rm d} \mathbb P + 2 \int_{\Omega_j'} N^{\frac4{p}} \, {\rm d} \mathbb P \\
\le \delta^2 N^{-2} + 2 N^{-1} \Ex |\xi_j|^p + 4N \mathbb P(\Omega_j')\,,
\end{multline*}
whence
\[
\sum_{j=1}^N \Ex |\xi_j - \eta_j'|^2 \le \delta^2 N^{-1} + 2N^{-1}\delta + 4N \cdot
2\delta N^{-2} \le 11 \delta N^{-1}\,,
\]
which together with $\Ex | \Theta |^2 \le \delta $ implies that the $L^2$-distance between
$\xi$ and $\sum_j \eta_j'$ does not exceed $ 5\sqrt{\delta} $. In particular, the variance
of the sum $\sum_j \eta_j'$ is close to the variance of $\xi$, that is, to $ 1 $.

Now, it remains to show that the independent random variables $\eta_j'$ satisfy the conditions
of Lyapunov's version of the central limit theorem. It remains to check that $\sum_j \Ex |\eta_j'|^p$
is small when $\delta$ is small. Indeed,
\[
\Ex |\eta_j'|^p = \int_{\Omega_j'} |\eta_j'|^p\, {\rm d}\mathbb P +
\int_{\Omega\setminus\Omega_j'} |\eta_j'|^p\, {\rm d}\mathbb P = J_1 + J_2\,.
\]
We have $J_1 \le N^2 \mathbb P(\Omega_j')$ because $|\eta_j'|^p \le N^2$. We also have
\[
J_2 \le \Ex \left[ |\xi_j| + \delta N^{-1} \right]^p
\le 2^p \left( \Ex |\xi_j|^p + \delta^p N^{-p}\right)\,.
\]
Thus, we finally get
\[
\Ex |\eta_j'|^p \le N^2 \mathbb P(\Omega_j') + 2^p \Ex |\xi_j|^p + 2^p \delta^p N^{-p}\,,
\]
and
\[
\sum_{j=1}^N \Ex |\eta_j'|^p \le N^2 \cdot 2\delta N^{-2} + 2^p \delta + 2^p \delta^p N^{1-p} \le 2(2^p+1)\delta\,,
\]
finishing the proof. \hfill $\Box$

\medskip
In our situation,
by Lemmas~\ref{lem_channels} and \ref{lem_FuzzySq}, the variance of the random variable
$\Theta$ is asymptotically negligible compared to  the variance of
$\overline{n}(R, h)$ as $R\to\infty$.
This gives us condition 1) in Lemma~\ref{lemma_prob}.

Now, we check that the sum of the $p$-th moments of the random variables $ \xi^{\mathcal Q} (R) $
is asymptotically negligible compared to the $\tfrac12 p$-th power of the variance of
$\overline{n}(R, h)$. This will give us condition 2) in Lemma~\ref{lemma_prob}. We choose
sufficiently large $p=p(\epsilon)$ (we shall have $p(\epsilon)\to \infty$ as $\epsilon\to 0$).

\subsection{The $p$-th moment estimate}
We have
\begin{multline*}
\Ex  \left| \xi_{\tt L}^{\mathcal Q} \right|^p =
\Ex \left| \frac1{2\pi R^2} \int_{R \mathcal Q} g\left( \frac{x}{R} \right) U(x)\, {\rm d} A(x)  \right|^p \\
\le \frac{\| g \|_\infty^p}{ (2\pi R)^{2p}} \, A(R\mathcal Q)^{p-1}
\int_{R \mathcal Q} \Ex \left| U(x) \right|^p \, {\rm d} A(x) \\
= \frac{\| g \|_\infty^p}{ (2\pi R)^{2p}} \, A(R\mathcal Q)^p \, \Ex \left| U(0) \right|^p
\end{multline*}
(in the last line, we used the distribution invariance of the random potential $U$).
Using estimate~\eqref{eq_g} and recalling that the area of the square $R \mathcal Q$
is $R^{2\gamma}$, we get
\[
\Ex  \left| \xi_{\tt L}^{\mathcal Q} \right|^p  \le C^p \| h \|_{C^\alpha}^p
\left[ \frac{R^{2-\alpha}}{R^2} R^{2\gamma} \right]^p = C(p, h) R^{(2\gamma-\alpha)p}\,.
\]

The $p$-th moment of the high-frequency component $\Ex  \left| \xi_{\tt H}^{\mathcal Q} \right|^p $
has the same upper bound. Indeed,
\begin{multline*}
\Ex  \left| \xi_{\tt H}^{\mathcal Q} \right|^p =
\Ex \left| \frac1{2\pi } \int_{\R^2}
h_{\tt H}^{\mathcal Q}\left( \frac{x}{R} \right) \Delta U(x)\, {\rm d} A(x)  \right|^p \\
\le \| h_{\tt H}^{\mathcal Q}\|_\infty^p
\Ex \left( n_F (R \widetilde{\mathcal Q}) + A(R \widetilde{\mathcal Q})\right)^p \\
\le C(p, h) R^{-p\alpha} \left(
\Ex n_F (R \widetilde{\mathcal Q})^p + A(R \widetilde{\mathcal Q})^p
\right)
\end{multline*}
where, as before, $\widetilde{\mathcal Q}$  denotes the support of the function
$\widetilde{\done}_{\mathcal Q}$.
We've used that the positive part of the signed measure $\tfrac1{2\pi} \Delta U \, {\rm d} A(x)$
is the counting measure $n_F$ of zeroes of $F$ while the negative part is the area measure
with the scaling coefficient $\tfrac12 $.

By distribution invariance of zeroes, estimating the $p$-th moment $\Ex n_F (R \widetilde{\mathcal Q})^p$,
we may assume that the square $\mathcal Q$ is centered at the origin. We take $\rho \simeq R^{\gamma}$
so that $\widetilde{\mathcal Q} \subset \{ |z| \le \rho \}$ and denote by
$n_F(\rho)$ the number of zeroes of $F$ in the disk $\{ |z|\le \rho \}$. Then
$\Ex n_F (R \widetilde{\mathcal Q})^p \le \Ex n_F(\rho)^p$.
By Jensen's formula,
\[
n_F (\rho) \le \frac1{2\pi} \int_0^{2\pi} \log^+ | F(e\rho e^{{\rm i}\theta}) |\, {\rm d}\theta
+ \log^-|F(0)|\,,
\]
and
\[
n_F(\rho)^p \le 2^p \left\{
\frac1{2\pi} \int_0^{2\pi} \bigl( \log^+ | F(e\rho e^{{\rm i}\theta}) | \bigr)^p \, {\rm d}\theta
+ \bigl( \log^- |F(0)| \bigr)^p \right\}\,.
\]
We have
\[
\log^+ | F(e\rho e^{{\rm i}\theta}) | \le \log^+ | F^*(e\rho e^{{\rm i}\theta}) | + \frac12 (e\rho)^2\,,
\]
whence,
\[
\bigl( \log^+ | F(e\rho e^{{\rm i}\theta}) | \bigr)^p
\le C(p) \left\{ \bigl( \log^+ | F^*(e\rho e^{{\rm i}\theta}) | )^p + \rho^{2p} \right\}\,.
\]
Since the distribution of the normalized function $F^*(z)$ does not depend on $z$,
and $F^*(0)=F(0)$ is a standard complex Gaussian random variable,
we have \[ \Ex \left\{ \bigl( \log^+ | F^*(e\rho e^{{\rm i}\theta}) | )^p \right\} =
\Ex \left\{ \bigl( \log^+ |F(0)| \bigr)^p\right\}. \]
Therefore, since $\rho>1$,
\[ \Ex n_F(\rho)^p \le C(p) \left[ \rho^{2p} + \Ex \left| \log |F(0)| \, \right| \right]^p
\le C_1(p) \rho^{2p}\,, \]
and recalling that $\rho \simeq R^{\gamma}$, we finally get
\[
\Ex  \left| \xi_{\tt H}^{\mathcal Q} \right|^p \le C(p, h) R^{p(2\gamma-\alpha)}\,.
\]

Since the total number of squares $\mathcal Q$ is of order $R^{2-2\gamma}$, we get
\[
\sum_{\mathcal Q} \left( \Ex \left| \xi_{\tt L}^{\mathcal Q} \right|^p
+ \Ex \left| \xi_{\tt H}^{\mathcal Q} \right|^p \right)
\le C(p, h) R^{2-2\gamma} \cdot R^{(2\gamma-\alpha)p}\,.
\]
On the other hand, by the assumption of the theorem we are proving,
\[
\sigma (R, h)^p \gtrsim c(p, h)R^{(\epsilon - \alpha)p}\,.
\]
Since $0<2\gamma< \epsilon$, choosing $p$ sufficiently big, we will have
\[
\sum_{\mathcal Q} \Ex \left| \xi^{\mathcal Q} \right|^p = o(1)\, \sigma (R, h)^p\,, \qquad R\to\infty\,.
\]

\subsection{Correcting the potential} Now, we combine the low and high frequency
terms. We set
\[
g^{\mathcal Q} = g \done_{\mathcal Q} + \Delta h_{\tt H}^{\mathcal Q}\,,
\]
and note that the functions $g^{\mathcal Q}$ have disjoint supports (provided that $R$
is big enough) and that
\[
\| g^{\mathcal Q} \|_\infty \le C(h) R^{M(\alpha)+2}
\]
with $M(\alpha)$ taken from \eqref{eq_prelim_2}. Then
\[
\xi^{\mathcal Q} (R) = \xi_{\tt L}^{\mathcal Q}(R) + \xi_{\tt H}^{\mathcal Q}(R)
=\frac1{2\pi R^2} \int_{\R^2} g^{\mathcal Q}\left( \frac{x}{R} \right) U(x)\, {\rm d} A(x)\,.
\]
It remains to approximate the random variables
$\xi^{\mathcal Q} (R)$ by independent ones. This will
be done in two steps. First, we correct the
random potential $U=\log|F^*|$ replacing the logarithmic kernel by a bounded
Lipschitz function $L$ with the Lipschitz norm $ R^{\frac12 B}$.

Given $R\ge 2$ and $B \ge 1$, let
\[
L(t) =
\begin{cases}
\ \ \log t, &{\rm when\ } |\log t| \le \frac12 B \log R, \\
\ \ \frac12 B \log R, &{\rm when\ } \log t > \frac12 B \log R, \\
-\frac12 B\log R, &{\rm when\ } \log t < -\frac12 B\log R,
\end{cases}
\]
let $U_{\tt cor} = L(F^*)$ be the corrected potential, and
let
\[
\xi_{\tt cor}^{\mathcal Q}(R) = \frac1{2\pi R^2} \int_{R \mathcal Q} g^{\mathcal Q}\left( \frac{x}{R} \right)
U_{\tt cor}(x)\, {\rm d} A(x)\,.
\]
be the corrected random variables.

\begin{lemma}\label{lemma_correction_error}
\[
\sum_{\mathcal Q} \Ex \left| \xi_{\tt cor}^{\mathcal Q}
- \xi^{\mathcal Q}  \right|
\lesssim R^{-12}\,,
\]
provided that $B > M(\alpha)+15$ and that $R$ is big enough.
\end{lemma}

\medskip\par\noindent{\em Proof of Lemma~\ref{lemma_correction_error}:}
We have
\begin{multline*}
\sum_{\mathcal Q} \left| \xi_{\tt cor}^{\mathcal Q}
- \xi^{\mathcal Q}  \right|
\le \frac1{2\pi R^2} \int_{\R^2} \left( \sum_{\mathcal Q}
\left| g^{\mathcal Q} \right| \right) \, \left| U_{\tt cor}(x) - U(x) \right|\, {\rm d} A(x) \\
\le \frac1{2\pi R^2}\,\sup_{\mathcal Q} \| g^{\mathcal Q}  \|_\infty
\, \int_{\{ |x|\le 2R \}} \left| U_{\tt cor}(x) - U(x) \right|\, {\rm d} A(x)\,.
\end{multline*}
We've used that the supports of the functions  $g^{\mathcal Q}$ are disjoint and
contained in the disk $\{ |x|\le 2R\}$.
The potentials $U_{\tt cor}(x)$ and $U(x)$ are different only when
$|F^*(x)|\le R^{-B/2}$ or $|F^*(x)|\ge R^{B/2}$. Hence, the right-hand side does not
exceed
\begin{multline*}
\frac{\sup_{\mathcal Q} \| g^{\mathcal Q}  \|_\infty}{2\pi R^2}\,
 \int_{\{ |x|\le 2R \}}  | U(x) |
\left( \done_{|F^*|\le R^{-B/2}} +  \done_{|F^*| \ge R^{B/2}} \right)
\, {\rm d} A(x) \\
\le C(h) R^{M(\alpha)} \int_{\{ |x|\le 2R \}} | U(x) |
\left( \done_{|F^*|\le R^{-B/2}} +  \done_{|F^*| \ge R^{B/2}} \right)
\, {\rm d} A(x) \,.
\end{multline*}
Here, we've used that $\| g^{\mathcal Q} \|_\infty
\le C(h) R^{M(\alpha)+2}$. Therefore,
\begin{multline*}
\sum_{\mathcal Q} \Ex \left| \xi_{\tt cor}^{\mathcal Q}
- \xi^{\mathcal Q}  \right| \le C(h) R^{M(\alpha)} \\
\cdot  \int_{\{ |x|\le 2R \}} \Ex \left\{ | U(x) |
\left( \done_{|F^*|\le R^{-B/2}} +  \done_{|F^*| \ge R^{B/2}} \right) \right\}
\, {\rm d} A(x) \,.
\end{multline*}

Due to the translation invariance of the distribution of $|F^*|$,
\begin{multline*}
\Ex \left\{ | U(x) |
\left( \done_{|F^*|\le R^{-B/2}} +  \done_{|F^*| \ge R^{B/2}} \right) \right\} \\
= \Ex \left\{ | \log |\zeta| |  \done_{|\zeta|\le R^{-B/2}} \right\}
+  \Ex \left\{  | \log |\zeta| |  \done_{|\zeta|\ge R^{B/2}} \right\}\,,
\end{multline*}
where $\zeta$ is a standard complex Gaussian random variable. The expectations on the right-hand
side are readily estimated:
\begin{multline*}
\Ex \left\{  | \log |\zeta| \, | \,\, \big| \,
|\zeta|\le R^{-B/2} \right\} \le \int_{\{|\zeta| \le R^{-B/2}\}} \log\frac1{|\zeta|}\, {\rm d} A(\zeta) \\
\lesssim \int_0^{R^{-B/2}} \bigl( \log\frac1{r} \bigr) r  {\rm d} r \lesssim B R^{-B} \log R\,,
\end{multline*}
and
\begin{multline*}
\Ex \left\{ | \log |\zeta| \, | \,\, \big| \,
|\zeta|\ge R^{B/2} \right\} \le \int_{\{|\zeta| \ge R^{B/2}\}} (\log|\zeta|) e^{-|\zeta|^2}\, {\rm d} A(\zeta) \\
\lesssim \int_{R^{B/2}}^\infty (\log r) e^{-r^2} r\, {\rm d} r \lesssim  e^{-R^B} B \log R\,.
\end{multline*}
Whence,
\[
\sum_{\mathcal Q} \Ex \left| \xi_{\tt cor}^{\mathcal Q}
- \xi^{\mathcal Q}  \right| \le B C(h) R^{M(\alpha)+2-B} \log R < R^{-12}\,,
\]
provided that $B> M(\alpha)+15$ and that $R$ is big enough.
\hfill $\Box$

\subsection{From the almost independent random variables to the independent ones}
At last, using Theorem~\ref{thm_almost_indep}, we approximate
the random variables $\xi_{\tt cor}^{\mathcal Q}$ by the independent ones $ \eta^{\mathcal Q}$.
We apply Theorem~\ref{thm_almost_indep}
to the compact sets $K = R \widetilde{\mathcal Q}$ and the values $\rho_{\mathcal Q} = R^{\frac12 \beta}$
which are much bigger than $\sqrt{\log \left(3+d(R\widetilde{\mathcal Q}) \right)}$. We assume that $R$ is
big enough. Then we get a collection of independent G.E.F.'s $F_{\mathcal Q}$ such that
for each square $\mathcal Q$,
\[
F^* = F^*_{\mathcal Q} + G^*_{\mathcal Q} \qquad {\rm on \ } R\widetilde{\mathcal Q}\,,
\]
with
\begin{equation}\label{eq_G_Q}
\mathbb P \left\{ \max_{R \widetilde{\mathcal Q}} |G^*_{\mathcal Q}| \ge e^{-R^{\beta}} \right\}
\lesssim \exp\bigl[ -e^{R^\beta} \bigr]\,.
\end{equation}

We set
\[
\eta^{\mathcal Q} (R) = \frac1{2\pi R^2} \int_{\R^2} g^{\mathcal Q}\left(\frac{x}{R}\right)
L( |F^*_{\mathcal Q}(x)| )\, {\rm d} A(x)\,,
\]
and check that
\[
\sum_{\mathcal Q} \mathbb P \left\{ \left| \xi^{\mathcal Q} - \eta^{\mathcal Q} \right|
\ge \delta N^{-2} \sigma \right\} \le \delta N^{-2}
\]
with $\delta = \tfrac1{R}$, $N \approx R^{2-2\gamma}$, and $\sigma (R, h) \ge c(h)R^{-1}$.
This will give us condition 3) of Lemma~\ref{lemma_prob}. With these values of the parameters
$\delta$, $N$, and $\sigma$, we have $ \delta N^{-2} \sigma > R^{-6} $
and $ \delta N^{-2} > R^{-5} $,
provided that $R$ is big enough. Hence, it suffices to check that
\[
\sum_{\mathcal Q} \mathbb P \left\{ \left| \xi^{\mathcal Q} - \xi_{\tt cor}^{\mathcal Q} \right|
\ge R^{-6} \right\} < R^{-5}
\]
and that
\[
\sum_{\mathcal Q} \mathbb P \left\{ \left| \xi_{\tt cor}^{\mathcal Q} - \eta^{\mathcal Q} \right|
\ge R^{-6} \right\} < R^{-5}\,.
\]
The first estimate follows from Lemma~\ref{lemma_correction_error}, so we need to check only
the second one.

We have
\begin{multline*}
\left| \xi_{\tt cor}^{\mathcal Q} - \eta^{\mathcal Q} \right|
\le
\frac1{2\pi R^2} \int_{\R^2} \bigl| g^{\mathcal Q}\left( \frac{x}{R}\right)\bigr|\,
\bigl| L(F^*_{\mathcal Q} (x) + G^*_{\mathcal Q} (x)) -
L(F^*_{\mathcal Q} (x)) \bigr| \, {\rm d} A(x) \\
\le \frac{\|g^{\mathcal Q} \|_\infty}{R^2} \cdot \| L \|_{\tt Lip}
\cdot \max_{R\widetilde{ \mathcal Q}} |G^*_{\mathcal Q}| \cdot A(R\mathcal Q) \\
\lesssim \frac{C(h) R^{M(\alpha)+2}}{R^2}
\cdot R^{\frac12 B} \cdot \max_{R\widetilde{\mathcal Q}} |G^*_{\mathcal Q}|
\cdot R^{2-2\gamma} \\
\stackrel{B>M(\alpha)+15}< R^{2B} \max_{R\widetilde{\mathcal Q}} |G^*_{\mathcal Q}|\,,
\end{multline*}
whence,
\begin{multline*}
\sum_{\mathcal Q } \mathbb P \bigl\{ \,  | \xi_{\tt cor}^{\mathcal Q} - \eta^{\mathcal Q}|
> R^{-6} \, \bigr\}
\le \sum_{\mathcal Q } \mathbb P \bigl(\, R^{2B} \max_{R \widetilde{\mathcal Q}}
|G_{\mathcal Q}^*| > R^{-B} \, \bigr) \\
\le \sum_{\mathcal Q } \mathbb P \bigl(\, \max_{R \widetilde{ \mathcal Q}}|G_{\mathcal Q}^*| > e^{-R^\beta} \,\bigr)
\lesssim N \exp\bigl[ -e^{R^\beta} \bigr],
\end{multline*}
which is much less than $R^{-5}$, provided that $R$
is big enough. This finishes off the proof of Theorem~\ref{thm_normality}.
\hfill $\Box$

\section{Abnormal test-functions}

By $\zeta$ and $\zeta_j$ we always denote standard complex Gaussian
random variables. Let $b=\Ex \log |\zeta|$. As above, we denote
$ \overline{n}(R, h) = n(R, h) - \Ex n(R, h) $.

\subsection{The function $h(x)=\log^-|x|$ is abnormal}
We start with a simple example of a function $h$ with abnormal fluctuations of linear statistics.
Though this function is unbounded, it may be regarded as a toy model for the main
example which we give in the next section.

In view of Jensen's integral formula, this
function is customary in the entire functions theory. We have
\[
n(R, h) = \sum_{a\in\cZ_F} \log^+ \frac{R}{|a|}
= \frac1{2\pi} \int_0^{2\pi} \log |F(Re^{{\rm i}\theta})|\, {\rm d}\theta - \log|F(0)|\,.
\]
Then
\[
\Ex n(R, h) = \frac1{2\pi} \int_0^{2\pi} \Ex \log |F(Re^{{\rm i}\theta})|\, {\rm d}\theta
- \Ex \log |F(0)|  = \Bigl( \frac12 R^2 + b \Bigr) - b = \frac12 R^2\,.
\]
Letting $F^*(z)=F(z)e^{-\frac12 |z|^2}$ and $\ell (\zeta) = \log|\zeta| - b$,
we see that
\begin{equation}\label{eq-Jensen}
\overline{n}(R, h) = \frac1{2\pi} \int_0^{2\pi} \ell \bigl( F^*(Re^{{\rm i}\theta}) \bigr)\,
{\rm d}\theta - \ell (F^*(0))\,.
\end{equation}
The variance of the first term on the right-hand side equals
\begin{multline*}
\frac1{4\pi^2} \int_0^{2\pi} \int_0^{2\pi}
\Ex \left\{ \ell\bigl( F^*(Re^{{\rm i}\theta_1}) \bigr) \ell\bigl( F^*(Re^{{\rm i}\theta_2} \bigr) \right\}
\, {\rm d}\theta_1 {\rm d}\theta_2 \\
\stackrel{{\rm Lemma}~\ref{lemma_covar}}= \frac1{16\pi^2} \sum_{\alpha\ge 1} \frac1{\alpha^2}
\int_0^{2\pi} \int_0^{2\pi}
e^{-\alpha R^2 |e^{{\rm i} \theta_1} - e^{{\rm i} \theta_2} |^2}\, {\rm d}\theta_1 {\rm d}\theta_2 \\
\le \frac1{8\pi} \sum_{\alpha\ge 1} \frac1{\alpha^2}
\int_0^{2\pi}
e^{-\alpha R^2 \sin^2\theta}\, {\rm d}\theta = O\left( \frac1{R}\right)\,, \qquad R\to\infty\,.
\end{multline*}
Therefore, the first term on the right-hand side of \eqref{eq-Jensen} can be
disregarded for large $R$'s.
Obviously, the fluctuations of the term $\ell (F^*(0))$ are not normal,
hence, the fluctuations of $n(R, h)$ are not normal as well. \hfill $\Box$

\subsection{Abnormal $C^\alpha$ functions with $0<\alpha<1$}
By $c_\alpha$ we denote various positive constants that depend
only on $\alpha$. We fix a radial function $\psi\in C^\infty_0(\R^2)$ that equals $1$ when $|x|\le 1$
and vanishes when $|x|\ge 2$, and put $h_\alpha (x) = |x|^\alpha \psi (x)$. Note
that for $|x|<1$, we have $\Delta h_\alpha (x) = c_\alpha |x|^{\alpha-2}$ with some $c_\alpha>0$.

We set $ \overline{U}(x) = \ell (F^*(x)) = \log |F^*(x)| - b $, and
introduce the random variable
\[
\xi_R = \int_{|x|\le R} |x|^{\alpha-2}\, \overline{U}(x)\, {\rm d} A(x)\,.
\]
By Remark~\ref{remark-var}, the second moment of the integral
$\displaystyle  \int_{A\le |x| \le B} |x|^{\alpha-2} \overline{U}(x) {\rm d} A(x)$ does not exceed
$\displaystyle C \int_A^B r^{2\alpha - 4} \cdot r\, {\rm d} r < c_\alpha A^{-2(1-\alpha)}$.
Therefore, for $R\to\infty$, the random variable $\xi_R$ converges in mean square to
the random variable
\[
\xi = \int_{\R^2} |x|^{\alpha-2} \overline{U}(x)\, {\rm d} A(x)\,.
\]

\begin{lemma}\label{lemma-nonconstant}
The random variable $\xi$ is not a constant one.
\end{lemma}

\medskip\par\noindent{\em Proof of Lemma~\ref{lemma-nonconstant}}:
Suppose that $\xi$ is a constant. Since $\Ex \xi = 0$, this means
that $\xi=0$ a.s.\,. Put
\[
\xi_R (y) = \int_{|x|\le R} |y-x|^{\alpha-2}\, \overline{U}(x)\, {\rm d} A(x)\,.
\]
By the translation invariance of the distribution of $\overline{U}$,
for each $y\in\R^2$,
\[
\lim_{R\to\infty} \Ex \bigl\{ \xi_R(y)^2 \bigr\} = 0\,.
\]
Since the random potential $\overline{U}(x)$ belongs to $L^2_{{\tt loc}}(\R^2)$,
the random functions $\xi_R(y)$ also belong to $L^2_{{\tt loc}}(\R^2)$.

We take a function $\chi\in C^\infty_0(\R^2)$ with positive Fourier transform,
and let $k = |x|^{\alpha-2} \ast \chi$.
Then
\[
\int_{|x|\le R}  k(x) \overline{U}(x)\, {\rm d}A(x)
= \int_{|y|\le R} \chi (y) \xi_R(y)\, {\rm d}A(y)  \stackrel{L^2}\to  0
\]
in mean square, when $R\to\infty$.
Thus,
\[
\int_{\R^2} k(x) \overline{U}(x)\, {\rm d} A(x) = 0\,,
\]
and by the translation invariance,
\[
\int_{\R^2} k(y-x) \overline{U}(x)\, {\rm d}A(x) \equiv 0\,, \qquad y\in\R^2.
\]
As before, the integrals converge in mean square.

The Fourier transform of the function $|x|^{\alpha-2}$ equals $c(\alpha)|\la|^{-\alpha}$
with $c(\alpha)>0$, so that $ \widehat{k}(\la) > 0 $.
Since $k$ is an $ L^2(\R^2) $-function with positive Fourier transform, any $L^2(\R^2)$-function
$ h(x) $ can be approximated in the $L^2(\R^2)$-norm by finite linear combinations of translations
$k(y_i-x)$ (this is a classical Wiener's theorem).
Then, by Remark~\ref{remark-var}, for {\em any} $L^2(\R^2)$-function $h$,
\[
\int_{\R^2} h(x) \overline{U}(x)\, {\rm d}A(x) = 0\,.
\]
In turn, this yields the absurd conclusion that $\overline{U}$ is the zero function. \hfill $\Box$

\begin{proposition}\label{thm-abnormal} \mbox{}

\smallskip\par\noindent {\rm (i)}
The random variables $R^\alpha \overline{n}(R, h_\alpha)$ converge in mean square
to $c_\alpha \xi$ as $R\to\infty$.

\smallskip\par\noindent {\rm (ii)}
The random variable $\xi$ is not a Gaussian one.

\smallskip
Hence, $R^\alpha \sigma (R, h_\alpha)$ converges to a positive limit, and
the fluctuations of the linear statistics $n(R, h_\alpha)$ are not asymptotically
normal.
\end{proposition}

\medskip\par\noindent{\em Proof of (i) in Proposition~\ref{thm-abnormal}}:
Since $\Delta h_\alpha$ is an $L^1$-function, we have
\[
\overline{n}(R, h_\alpha) = \frac1{2\pi R^2} \int_{\R^2} (\Delta h_\alpha)\bigl(\tfrac{x}{R}\bigr)
\overline{U}(x)\, {\rm d} A(x)\,.
\]
Therefore,
\[
R^\alpha \overline{n}(R, h_\alpha) = c_\alpha \xi_R +
\frac1{2\pi R^{2-\alpha}} \int_{R \le |x| \le 2R} (\Delta h_\alpha)\bigl(\tfrac{x}{R}\bigr)
\overline{U}(x)\, {\rm d} A(x)\,.
\]
Using Remark~\ref{remark-var}, we see that the variance of the second term
on the right-hand side does not exceed
\[
c_\alpha R^{2\alpha-4} \int_{R\le |x| \le 2R}
\left| \Delta h_\alpha \bigl( \tfrac{x}{R}\bigr) \right|^2\, {\rm d} A(x)
\le c_\alpha R^{2\alpha-2} \to 0
\]
for $R\to\infty$.
Hence, $R^\alpha \overline{n}(R, h_\alpha)$ converges in mean square to $c_\alpha \xi$.
\hfill $\Box$

\medskip

In what follows, we use some elementary relations:
\begin{lemma}\label{lemma-5.2.elementary} For $t>0$, we have
\begin{equation}\label{eq-5.2-elementary1}
\Ex \left\{ |\zeta|^t \right\} = \Gamma \bigl( \tfrac12 t + 1 \bigr)\,,
\end{equation}
\begin{equation}\label{eq-5.2-elementary2}
\Gamma \bigl( \tfrac12 t + 1 \bigr) \le e^{bt+Ct^2}\,,
\end{equation}
and
\begin{equation}\label{eq-5.2-elementary3}
e^{-bt} \Ex \left\{ |\zeta|^t \right\} \le e^{Ct^2}\,,
\end{equation}
where $C$ is a sufficiently large positive
numerical constant.
\end{lemma}

\par\noindent{\em Proof of Lemma~\ref{lemma-5.2.elementary}:}
For $t \ge 0$, we have
\[
\Ex \left\{ |\zeta|^t \right\} = \frac1{\pi} \int_{\C} |z|^t e^{-|z|^2}\, {\rm d} A(z) =
\int_0^\infty r^{t/2} e^{-r}\, {\rm d} r = \Gamma \left( \tfrac{t}2 + 1\right).
\]
This gives us~\eqref{eq-5.2-elementary1}. Then
it is easy to check that for $s\ge 0$,
\[
\log \Gamma (1+s) \le \frac{\Gamma'(1)}{\Gamma (1)}\, s + C s^2
\]
with a sufficiently big positive numerical constant $C$.
It remains to note that
\[
b = \Ex \left\{ \log |\zeta| \right\} = \frac1{\pi} \int_\C \left( \log |z| \right)
\, e^{-|z|^2}\, {\rm d}A(z)
= \frac12 \int_0^\infty \left( \log s \right) \, e^{-s}\, {\rm d}s =
\frac12 \frac{\Gamma'(1)}{\Gamma (1)}\,,
\]
completing the proof of \eqref{eq-5.2-elementary2}. \hfill $\Box$

\medskip

The idea behind the proof of (ii) in Proposition~\ref{thm-abnormal} is also simple.
At each point $x\in \R^2$, the random potential $\U(x)$ is distributed like
$\log |\zeta| - b$, whence by \eqref{eq-5.2-elementary1}
$ \log \Ex \left\{ e^{t \U(x)} \right\} \le Ct\log t $ for $t \to + \infty$.
Since the random variable $\xi$ is a weighted average of the random potential
$\overline{U} (x)$, we expect that $ \Ex \left\{ e^{t\xi} \right\} \le e^{o(t^2)} $ for $t\to +\infty$,
which forbids the random variable $\xi$ to be a Gaussian one.
To implement this programme, we need to estimate some Laplace transforms.

\begin{lemma}\label{lemma-5.2.charfunct}
Suppose $\phi$ is a non-negative function on $\R^2$. Then

\smallskip\par\noindent{\rm (i)}
\[
\Ex \exp \left[ \int_{\R^2} \phi \U  \right] \le e^{-b \|\phi\|_1} \Gamma \bigl( \tfrac12 \|\phi\|_1 +1 \bigr);
\]

\smallskip\par\noindent{\rm (ii)} for each $L>0$,
\[
\Ex \exp \left[ \int_{\R^2} \phi \U  \right] \le e^{\frac12 \delta (L) \|\phi\|_1 + CL^2 \| \phi\|_2^2}\,,
\]
where $\delta(L) = \sum_{k\in\Z^2\setminus\{0\}} e^{-L^2|k|^2/2}$.
\end{lemma}

\par\noindent{\em Proof of (i) in Lemma~\ref{lemma-5.2.charfunct}:}
By Jensen's inequality,
\[
\exp \left[ \int_{\R^2} \phi \U\, {\rm d} A \right]
\le \int_{\R^2} e^{\| \phi \|_1 \U}\, \frac{\phi\, {\rm d}A}{\|\phi\|_1}\,.
\]
Hence,
\begin{multline*}
\Ex \exp \left[ \int_{\R^2} \phi \U \, {\rm d}A \right]
\le \int_{\R^2} \Ex \left\{ e^{ \| \phi \|_1 (\log |F^*|-b) } \right\}\, \frac{\phi {\rm d}A}{\|\phi\|_1}
\\
= e^{-b\|\phi\|_1} \Ex \left\{ |\zeta|^{\| \phi \|_1} \right\}
= e^{-b \|\phi\|_1} \Gamma \bigl( \tfrac12 \|\phi\|_1 +1 \bigr)\,,
\end{multline*}
proving (i). \hfill $\Box$

\medskip The second estimate in Lemma~\ref{lemma-5.2.charfunct}
is more delicate. It's proof uses two lemmas. The first lemma is
classical (see Theorem~26 in~\cite{HLP} where it is attributed to
Ingham and Jessen).

\begin{lemma}\label{lemma-5.2-conv}
Suppose $X$ is a measure space with a probability measure $\mu$, and $f\colon X\to\R$
is a random function on $X$. Then
\begin{equation}\label{eq-5.2-convex}
\log\, \Ex\, \exp \left[ \int_X f\, {\rm d}\mu \right] \le \int_X \log\, \Ex \left\{  \exp f \right\}\,
{\rm d}\mu\,.
\end{equation}
\end{lemma}

\medskip\par\noindent{\em Proof of Lemma~\ref{lemma-5.2-conv}:} Consider the functional
$\theta \mapsto \log\, \Ex\, e^\theta$ on real-valued random variables. By H\"older's inequality,
it is convex: for $0\le t\le 1$, we have
\begin{multline*}
\log\, \Ex \bigl\{  e^{ t \theta_1 + (1-t) \theta_2)} \bigr\}
\le \log\, \left\{ \left( \Ex \left\{ e^{\theta_1} \right\} \right)^t
\cdot \left( \Ex \left\{ e^{\theta_2} \right\} \right)^{1-t} \right\} \\
= t \log \bigl\{ \Ex\, e^{\theta_1} \bigr\} +
(1-t) \log\, \Ex \bigl\{ e^{\theta_2} \bigr\}.
\end{multline*}
Then \eqref{eq-5.2-convex} follows from  Jensen's inequality. \hfill $\Box$

\begin{lemma}\label{lemma-5.2.4}
Let $\zeta_i$ be standard complex Gaussian random variables such that for every $i$, we have
$$
\sum_{j\colon j\ne i} \left| \Ex \left\{ \zeta_i \bar\zeta_j \right\} \right| \le\delta\,.
$$
Then for all $t_i\ge 0$, we have
\begin{equation}\label{eq-5.2.a}
\log \Ex \prod_i|\zeta_i|^{t_i}\le \frac\delta 2 \sum_i t_i + b \sum_i t_i
+ C\sum_i t_i^2\,.
\end{equation}
\end{lemma}

\medskip\par\noindent{\em Proof of Lemma~\ref{lemma-5.2.4}:}
Since the matrix with the entries
$$
a_{ij}=\begin{cases}
\delta, & i=j\\
-\Ex \left\{ \zeta_i\bar\zeta_j \right\}, & i\ne j
\end{cases}
$$
is non-negative definite, we can find complex Gaussian random variables $\eta_i$
independent of all $\zeta_i$ such that $\Ex \left\{ \eta_i\bar\eta_j \right\} = a_{ij}$. For $z\in\C$,
put $\la_i(z) = \zeta_i +z \eta_i$. When $z\in\T$, the random variables
$\frac{\la_i(z)}{\sqrt{1+\delta}}$ are independent standard complex Gaussians,
therefore
\[
\Ex \prod_i|\la_i(z)|^{t_i} = \prod_i \Ex |\la_i(z)|^{t_i}\,.
\]
By Lemma~\ref{lemma-5.2.elementary},
\[
\Ex |\la_i(z)|^{t_i} \le (1+\delta)^{\frac12 t_i}\, e^{ b t_i + Ct_i^2}\,.
\]
Thus
\[
\log \Ex \prod_i|\la_i(z)|^{t_i} \le \frac\delta 2 \sum_i t_i + b \sum_i t_i
+ C\sum_i t_i^2\,.
\]

On the other hand, $\prod_i|\la_i(z)|^{t_i}$ is a subharmonic function of $z$, so
$$
\prod_i|\zeta_i|^{t_i}=\prod_i|\la_i(0)|^{t_i}\le\int_{\T}\prod_i|\la_i(z)|^{t_i}\,dm(z)
$$
where $m$ is the Haar measure on $\T$. Taking the expectation of both sides,
we complete the proof. \hfill $\Box$

\medskip Now we can complete the proof of Lemma~\ref{lemma-5.2.charfunct}.

\medskip\par\noindent{\em Proof of (ii) in Lemma~\ref{lemma-5.2.charfunct}:} We have
\begin{multline*}
\Ex \exp \left[ \int_{\R^2} \phi \U  \right] =
\Ex \exp \left[ \frac1{L^2} \int_{[0, L]^2} L^2 \sum_{k\in\Z^2 }\phi (x+kL)  \U (x+kL)  \right] \\
\stackrel{\rm\eqref{eq-5.2-convex}}\le
\exp \left[ \frac1{L^2} \int_{[0, L]^2} \log\, \Ex\, e^{L^2 \sum_k \phi (x+kL) \U (x+kL)}
\right] \\
= e^{-b \| \phi \|_1}
\exp \left[ \frac1{L^2} \int_{[0, L]^2}
\log\, \Ex \prod_{k\in\Z^2} |F^*(x+kL)|^{L^2 \phi (x+kL)} \right].
\end{multline*}
To estimate the expectation of the product on the right-hand side, we apply
Lemma~\ref{lemma-5.2.4}. We have
\[
\sum_{k'\colon k'\ne k} \left| \Ex \left\{ F^*(x+kL) \overline{F^*(x+k'L)} \right\} \right|
= \sum_{k'\colon k'\ne k} e^{-L^2 |k'-k|^2/2}
= \sum_{k\in \Z^2\setminus\{0\}} e^{-L^2 |k|^2/2}\,.
\]
Denoting the sum on the right-hand side by $\delta (L)$, and applying Lemma~\ref{lemma-5.2.4},
we get
\begin{multline*}
\log\, \Ex \prod_{k\in\Z^2} |F^*(x+kL)|^{L^2 \phi (x+kL)}
\\
\le \frac12 \delta (L) L^2 \sum_{k\in\Z^2} \phi (x+kL) + CL^4 \sum_{k\in\Z^2} \phi^2(x+kL)
+ bL^2 \sum_{k\in\Z^2} \phi (x+kL)\,,
\end{multline*}
whence
\begin{multline*}
\log \Ex \exp \left[ \int_{\R^2} \phi \U  \right] \le -b\|\phi\|_1 \qquad \qquad \\
\qquad + \frac1{L^2} \int_{[0, L]^2}
\Bigl(
\frac12 \delta (L) L^2 \sum_{k\in\Z^2} \phi (x+kL) + CL^4 \sum_{k\in\Z^2} \phi^2(x+kL)
+ bL^2 \sum_{k\in\Z^2} \phi (x+kL)
\Bigr) {\rm d}x \\
= \frac12 \delta (L) \| \phi\|_1 + CL^2 \| \phi \|_2^2\,,
\end{multline*}
completing the proof of Lemma~\ref{lemma-5.2.4}.
\hfill $\Box$

\medskip To prove that the random variable $\xi$ is not Gaussian,
we use a simple probabilistic lemma:

\begin{lemma}\label{lemma-5.2.3}
Let $\eps_j>0$ satisfy $S=\sum_{j=0}^\infty\eps_j<+\infty$.
Suppose that real random variables $\xi_j$ satisfy
\begin{equation}\label{eq-cond1}
\log \Ex e^{t\xi_j} \le o(t^2) \qquad \text{for all } j\ge 0, \text{and } t\to\infty\,,
\end{equation}
and
\begin{equation}\label{eq-cond2}
\log \Ex e^{t\xi_j}\le \eps_j^2 t^2
\qquad \text{for all } j\ge 1, \text{and all } t\ge 1.
\end{equation}
Then
$$
\log \Ex e^{t\sum_{j\ge 0}\xi_j} = o(t^2) \quad\text{as }t\to+\infty\,.
$$
In particular, this implies that $\sum_{j\ge 0}\xi_j$ is not normal.
\end{lemma}

\medskip\par\noindent{\em Proof of Lemma~\ref{lemma-5.2.3}:}
Let $\delta_j=\eps_j/S$. By H\"older's inequality, we have
$$
\Ex e^{t\sum_{j\ge 0}\xi_j } \le
\prod_{j\ge 0} \left\{ \Ex e^{ \delta_j^{-1} t \xi_j } \right\}^{\delta_j}
=\prod_{0\le j\le N}\cdot\prod_{j>N}\, \le e^{o(t^2)}\cdot e^{t^2 S \sum_{j>N}\eps_j}\,.
$$
It remains to note that $\sum_{j>N}\eps_j\to 0$ as $N\to+\infty$.
\hfill $\Box$

\medskip At last, everything is ready to show that the random variable $\xi$ is
not a Gaussian one, and thus to complete the proof of Proposition~\ref{thm-abnormal}.

\medskip\par\noindent{\em Proof of (ii) in Proposition~\ref{thm-abnormal}:}
We split the integral $\xi$ into the sum of the
integrals
\[
\xi_0=\int_{|x|\le 1} |x|^{\alpha-2}\U(x)\,{\rm d}A(x),
\quad {\rm and} \quad
\xi_j=\int_{2^{j-1}<|x|\le 2^j} |x|^{\alpha-2}\U(x)\, {\rm d}A(x) \ {\rm for} \  j\ge 1\,.
\]
To show that $\xi$ is not a Gaussian random variable, we check that
$\xi_j$'s satisfy assumptions~\eqref{eq-cond1} and \eqref{eq-cond2}
of Lemma~\ref{lemma-5.2.3}.
By estimate (i) in Lemma~\ref{lemma-5.2.charfunct}, condition~\eqref{eq-cond1} holds
for the random variables $\xi_j$, $j\ge 0$, even with $O(t\log t)$ on the right-hand
side.

To check condition~\eqref{eq-cond2} for $\xi_j$'s with $j\ge 1$, we apply estimate
(ii) in Lemma~\ref{lemma-5.2.charfunct} with
$ \phi_j(x) = t |x|^{\alpha-2} \done_{\{ 2^{j-1}<|x|\le 2^j \}} (x) $. Note that
\[
\| \phi_j \|_1 \le C_\alpha t 2^{\alpha j} \quad {\rm and} \quad   \| \phi_j \|_2^2 \le C_\alpha
t^2 2^{-2(1-\alpha)j}\,.
\]
Also note that for $L\ge 1$,  $\delta (L) = \sum_{k\in\Z^2\setminus\{0\}} e^{-L^2|k|^2/2} \le C e^{-cL^2}$
with positive numerical constants $c$ and $C$. Thus, estimate (ii) in Lemma~\ref{lemma-5.2.charfunct}
gives us the upper bound
\[
\log\, \Ex\, e^{t\xi_j} \le C_\alpha t^2 \left( e^{-cL^2} 2^{\alpha j} + L^2 2^{-2(1-\alpha)j} \right)
\]
valid for $t\ge 1$ and $L\ge 1$. Taking there $L=j$, we see that
$ \log \Ex e^{t\xi_j} \le C_\gamma t^2 2^{-2\gamma j} $ for $ t\ge 1 $
with any $0<\gamma<1-\alpha$. Thus, condition~\eqref{eq-cond2} is satisfied with
$\epsilon_j = C_\gamma 2^{-\gamma j}$, and the conclusion follows.
\hfill $\Box$

\end{document}